%% file: main.tex
\documentclass[sn-chicago,iicol]{sn-jnl}

\input{packages.tex}

\usepackage{hyperref}       
\usepackage{url}            
\usepackage{nicefrac}       
\usepackage{microtype}      
\usepackage{cleveref}       
\usepackage{doi}
\usepackage{stanli}

\usepackage{graphicx}%
\usepackage{stfloats}%
\usepackage{multirow}%
\usepackage{amsmath,amssymb,amsfonts}%
\usepackage{amsthm}%
\usepackage{mathrsfs}%
\usepackage[title]{appendix}%
\usepackage{xcolor}%
\usepackage{textcomp}%
\usepackage{manyfoot}%
\usepackage{booktabs}%
\usepackage{algorithm}%
\usepackage{algorithmicx}%
\usepackage{algpseudocode}%
\usepackage{listings}%

\begin{document}

\title[Topology optimization]{
	A nonlocal approach to graded surface modeling in topology optimization
}
\author*[1,2]{\fnm{Sukhminder} \sur{Singh}}\email{sukhminder.singh@fau.de}

\author[2]{\fnm{Lukas} \sur{Pflug}}\email{lukas.pflug@fau.de}

\author[1]{\fnm{Fabian} \sur{Wein}}\email{fabian.wein@fau.de}

\author[1]{\fnm{Michael} \sur{Stingl}}\email{michael.stingl@fau.de}

\affil*[1]{
	\orgdiv{Chair of Applied Mathematics (Continuous Optimization)},
	\orgname{Friedrich-Alexander-Universität Erlangen-Nürnberg},
	\orgaddress{
		\country{Germany}}%
}
\affil[2]{
	\orgdiv{Competence Unit for Scientific Computing (CSC)}, \par
	\orgname{Friedrich-Alexander-Universität Erlangen-Nürnberg},
	\orgaddress{
		\country{Germany}}%
}

\abstract{
	Additively manufactured structures often exhibit a correlation between their mechanical properties, such as stiffness, strength, and porosity, and their wall thickness.
	This correlation stems from the interplay between the manufacturing process and the properties of the filler material.
	In this study, we investigate the thickness-dependent effect on structural stiffness and propose a nonlocal integral model that introduces surface grading of Young's modulus to capture this phenomenon.
	We incorporate this model into topology optimization for designing structures with optimized compliance subject to a volume constraint.
	Notably, elastically degraded surfaces penalize excessively thin features, effectively eliminating them from the optimized design.
	We showcase the efficacy of our proposed framework by optimizing the design of a two-dimensional cantilever beam and a bridge.
}


\maketitle

\section{Introduction}

This paper introduces a novel approach to topology optimization, focusing on incorporating surface grading within the optimized structures.
Surface grading, a distinctive characteristic of additively manufactured structures, is intricately linked to both the manufacturing process and the mechanical properties inherent in the base material.
Notably, in powder-bed additive manufacturing techniques such as Electron Beam Melting (EBM) or Selective Laser Sintering (SLS), various material properties, including elastic modulus, Poisson's ratio, ultimate tensile strength, and porosity, exhibit non-uniformity across the thickness of the fabricated features.
This non-uniform distribution leads to thickness-dependent effective material characteristics, exerting a pronounced influence on the overall properties of the structure.
Such effects are particularly significant in thin-walled structures, where the thickness of the walls is comparable to the size of the laser beam spot used in the additive manufacturing process.

Multiple studies have investigated the thickness dependency of the material properties in additive manufacturing.
\citet{ThicknessDepenAlgard2016} studied the thickness dependency of microstructures in thin-walled titanium parts manufactured by EBM, revealing a significant impact on mechanical properties attributed to high surface roughness.
\citet{ThicknessDepenTasch2018} determined that structures possessing a thickness less than \SI{1}{\milli\meter} exhibited notable reductions in stiffness, ultimate tensile strength, and elongation at break.
\citet{ThicknessDepenSindin2020} demonstrated evidence of thickness dependency of the material properties of laser-sintered polyamide~$12$ (PA$12$) specimens, with the degree of variation in mechanical properties changing as wall thickness decreases.
Subsequently, the authors developed finite element models for the thickness dependent Young's modulus and Possion's ratio of shell structures, calibrated from experiments with laser-sintered short-fiber-reinforced PA12 material
~\citep{ConsideringInhSindin2021,MaterialModellSindin2021,StructuralRespSindin2021}.
The thickness dependency was modeled using polynomial functions of thickness for Young's modulus and Poisson's ratio.
On the other hand, \citet{jaksch2022thin} introduced a nonlocal integral model to incorporate surface grading of Young's modulus, simulating similar thickness-dependent effect.


Density-based topology optimization has emerged as a robust and versatile technique for designing structures with optimal material distribution. Conventional SIMP (Solid Isotropic Material with Penalization) approaches~\citep{TopologyOptimiBendso2004} penalize intermediate densities to promote a near-binary material distribution.
Furthermore, a Heaviside projection filter~\citep{AchievingMinimGuest2004} is incorporated to introduce a characteristic length scale in the design space, enabling feature size control and mitigating mesh dependency issues.

Recently, many \emph{coating} filters have been introduced in the literature to model surface layer effects in optimized structures.
\citet{TopologyOptimiClause2015} presented topology optimization framework of structures with stiff coating.
Subsequent to this investigation, various adaptations of the filter were explored in the literature \citep{ANewCoatingFYoon2019, HomogenizationGroen2019,AProjectionBaLuoY2019,TopologyOptimiSuresh2020}.
\citet{ALevelSetMetWang2018} introduced a level-set method to model coated structures within shape and topology optimization framework.

Additionally, \citet{TopologyOptimiTuna2022} conducted topology optimization of two-dimensional plates, incorporating Eringen's nonlocal theory of elasticity~\citep{TheoryOfNonloEringe1987}.
However, the distinctive aspect of our study lies in the departure from employing the nonlocal theory to calculate stress at individual points.
Instead, we employ an integral formulation for the material constant, enabling the computation of a spatially varying elastic modulus.

This study aims to develop a density-based topology optimization framework, particularly suited to thin-walled structures, considering the thickness-dependent nature of material properties.
The nonlocal integral model, as described in \cite{jaksch2022thin}, is utilized to capture the thickness dependence.
This model effectively considers surface layer effects and seamlessly integrates into topology optimization algorithm.

Our formulation, as explained in the following sections, penalizes slender structural features while favoring more robust and thicker counterparts in the resulting optimized design.
To avoid unwanted effects such as checkerboard pattern and mesh dependency in the results, and to enable control over the minimum feature size, we employ the conventional Heaviside projection filter.
The significant benefit of this approach is that it allows for the improved design of structures suitable for additive manufacturing processes, leading to a better alignment of material properties that depend on the variable thickness of parts.

The framework is particularly well-suited for topology optimization of lattice microstructures \citep{TopologyOptimiWuJu2021}, wherein the periodic cell dimensions typically fall within the range of \SI{5}{\milli\meter} to \SI{20}{\milli\meter}.
Conventional black-and-white designs may induce excessively thin features, while employing feature size control methods to mitigate this issue could compromise design flexibility.
To address these limitations, the current method aims to integrate actual material properties into the optimization process, enabling the generation of optimized designs devoid of artificial regularization.
This approach is expected to produce structures with a more natural and robust topology, better suited for practical applications.

The remainder of this paper is organized as follows.
Section~\ref{sec:problem_formulation} presents the problem formulation.
This includes the definition of the nonlocal material model and the optimization problem.
Section~\ref{sec:results} presents the results of the numerical experiments.
Finally, Section~\ref{sec:conclusion} concludes the paper.

\section{Problem formulation}
\label{sec:problem_formulation}

This section introduces the nonlocal approach developed by \citet{jaksch2022thin}, designed to simulate a variable Young's modulus across the specimen's thickness.

\subsection{Material model}
\label{sec:material_model}

\begin{figure}[t]
	\centering
	\begin{tikzpicture}
		\point{a}{0}{0};
		\point{c}{0}{0.5};
		\point{b}{5}{1};
		\draw[thick, fill=gray!10] (a) rectangle (b);
		\dimensioning{2}{a}{b}{3}[$t$]
		\draw[-latex] (1, 0.5) node[left] {$0$} -- (1, 1.5) node[above left] {$x$};
		\draw (0.9, 0.5) -- (1.1, 0.5);
		\foreach \y in {0,0.2,...,1} {
				\draw[thick, -latex] (5, \y) -- (5.5, \y);
				\draw[thick, -latex] (0, \y) -- (-0.5, \y);
			}
		\node at (3.7, 0.5) {$\Omega$};
	\end{tikzpicture}
	\caption{
		A one-dimensional specimen of thickness $t$ under tension.
	}
	\label{fig:elastic_bar}
\end{figure}
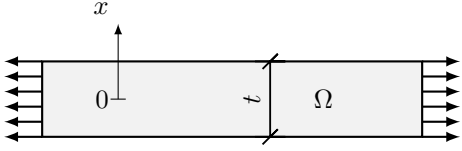
Consider $\Omega \in \mathbb{R}^d, d\in\{2, 3\}$ as a homogeneous, bounded physical domain.
Let $I_\Omega$ be an indicator function defining the interior of the domain, given by
\begin{equation}
	{I}_\Omega (\boldsymbol{x})
	\coloneqq
	\begin{cases}
		1 & \text{if } \boldsymbol{x} \in \Omega, \\
		0 & \text{otherwise}.
	\end{cases}
	\label{eq:indicator_function}
\end{equation}
Given a reference Young's modulus $E_0$ and
a kernel function $\psi: \mathbb{R}^d \rightarrow \mathbb{R}$, the Young's modulus at a point $\boldsymbol{x} \in \Omega$ is defined as
\begin{equation}
	E(\boldsymbol{x})
	\coloneqq
	\left[ \zeta + (1-\zeta) \phi(\boldsymbol{x}) \right] E_0,
\end{equation}
where $\phi(\boldsymbol{x})$ represents the nonlocal term, expressed as the convolution integral:
\begin{equation}
	\phi(\boldsymbol{x})
	\coloneqq
	\int_{\mathbb{R}^d} \psi(\boldsymbol{y}; \boldsymbol{x}) {I}_\Omega (\boldsymbol{y}) \, \mathrm{d}\boldsymbol{y}.
	\label{eq:nonlocal_youngs_modulus}
\end{equation}
This nonlocal term introduces a dependency on the spatial distribution of Young's modulus.
The parameter $\zeta \in [0, 1]$ serves as a fraction coefficient, offering control over the surface grading profile.
Particularly, for $\zeta=1$, the material exhibits a uniform Young's modulus distribution.
The kernel function $\psi$ captures the nonlocal effect, showcasing a decaying nature that vanishes beyond a certain limit.
It satisfies the normalization condition:
\begin{align}
	\int_{\mathbb{R}^d} \psi(\boldsymbol{y}; \boldsymbol{x})\, \mathrm{d}\boldsymbol{y} = 1 \quad \forall \boldsymbol{x} \in \Omega.
\end{align}
In this study, we adopt the following expression for the kernel function:
\begin{align}
	\psi(\boldsymbol{y}; \boldsymbol{x}) \coloneqq \frac{\max\left\{ \delta - \Vert \boldsymbol{y} - \boldsymbol{x} \Vert, 0 \right\}}{ \int_{\mathbb{R}^d } \max\left\{ \delta - \Vert \boldsymbol{z} - \boldsymbol{x} \Vert, 0 \right\} \, \mathrm{d}\boldsymbol{z} },
	\label{eq:kernel_function}
\end{align}
which is widely used in density-based topology optimization literature.
Here, $\delta > 0$ defines the region of influence of the kernel, determining the extent to which nonlocal interactions impact the material properties.
\begin{figure}[htbp]
	\centering
	\includegraphics[width=0.48\textwidth]{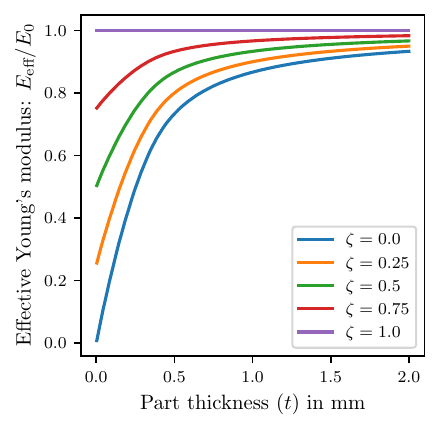}
	\caption{
		The effective Young's modulus (Eq.~\eqref{eq:effective_youngs_modulus_1d}) for one-dimensional elastic bar, with $\delta=\SI{0.4}{\milli\meter}$.
		The effective Young's modulus significantly degrades for thickness below \SI{1}{\milli\meter}.
	}
	\label{fig:effective_youngs_modulus_1d}
\end{figure}
\begin{figure*}[htbp]
	\begin{center}
		\includegraphics[width=\textwidth]{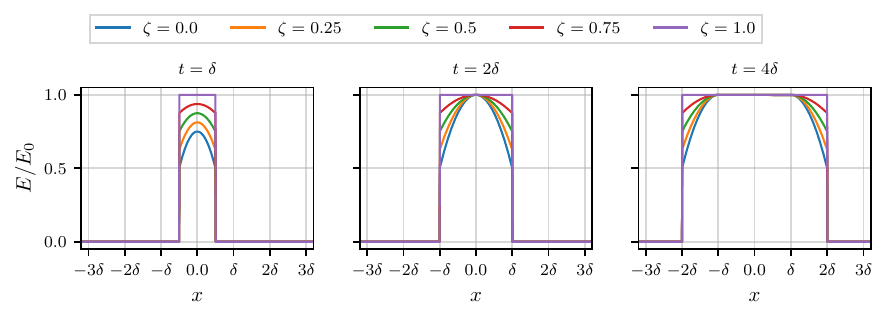}
	\end{center}
	\caption{
		Spatial distribution of Young's modulus across the thickness of the one-dimensional specimen, as illustrated in Figure~\ref{fig:elastic_bar}, showcasing variations corresponding to different specimen thickness values $t$.
	}
	\label{fig:youngs_modulus_profiles_1d}
\end{figure*}

To elucidate the impact of the selected material model on the effective Young's modulus concerning wall thickness, we analyze a one-dimensional specimen with thickness denoted by $t$, as illustrated in Figure~\ref{fig:elastic_bar}.
The homogenized effective Young's modulus at a point $x \in [-t/2, t/2]$, denoted as $E_\text{eff}$, is computed as the average of the varying Young's modulus $E$ across the thickness of the specimen, expressed as
\begin{align}
	E_\text{eff}(t) = \frac{1}{t} \int_{-t/2}^{t/2} E(x) \, \mathrm{d}x.
	\label{eq:effective_youngs_modulus_1d}
\end{align}
Figure~\ref{fig:effective_youngs_modulus_1d} shows the dependency of the effective Young's modulus on the thickness of the specimen, considering $\delta=\SI{0.4}{\milli\meter}$ \citep{jaksch2022thin} and varying values of the fraction coefficient $\zeta$.
The integrals in Eq.~\eqref{eq:nonlocal_youngs_modulus}, Eq.~\eqref{eq:kernel_function} and Eq.~\eqref{eq:effective_youngs_modulus_1d} are computed numerically using the midpoint method.
For values of $\zeta < 1$, the effective Young's modulus significantly decreases for thickness below \SI{1}{\milli\meter}.
For $\zeta=1$, the effective Young's modulus is independent of the specimen's thickness.

Figure~\ref{fig:youngs_modulus_profiles_1d} depicts the spatial variation of Young's modulus across the thickness of the one-dimensional specimen, highlighting distinct profiles associated with varying thickness values.
For specimens with a thickness $t \geq 2\delta$, the graded surface layer maintains a constant thickness of $\delta$.
In this scenario, the Young's modulus at the surface is $0.5E_0$
and increases gradually until reaching its maximum value of $E_0$ within the graded layer.
However, for specimen thicknesses below $2\delta$, the graded surface layer encompasses the entire thickness, leading to a decrease in both the maximum and minimum values of Young's modulus with diminishing thickness.

\subsection{Topology optimization}

The objective of this study is to design structures with graded surfaces in terms of elastic modulus, while simultaneously maintaining a near-binary material distribution.
For this purpose, the design domain is first discretized using $N$ finite elements,
followed by the application of a density-based topology optimization method within the spatially discrete framework.

In conventional topology optimization using the SIMP method, a piecewise-constant, scalar density field, represented by $\rho_e, e=1,\dots,N$, is employed.
Here, $\rho_e=0$ denotes a void element, $\rho_e=1$ indicates a solid element, and intermediate values of $\rho_e \in (0, 1) $ represent gray regions.
These intermediate regions are penalized in the optimization process (explained later), promoting a near-binary material distribution.
These values are termed ``densities'' as they are combined with elemental volumes $v_e, e=1,\dots,N$, to give the total mass of the structure.

However, a drawback of this penalization approach is the emergence of checkerboard pattern in the optimized design when piecewise-linear finite element interpolation is used for the displacement field.
To overcome this issue, a pseudo density field, given by $\mu_e, e=1,\dots,N$, is introduced.
The pseudo densities act as control variables in the optimization process, and lack a direct physical interpretation.
A convolution filter is applied to the pseudo density field which effectively mitigates the occurrence of checkerboard pattern~\citep{TopologyOptimiBruns1998}.

The introduction of the filter also establishes a length scale in the optimized design, addressing the issue of the optimizer converging to designs dependent on the resolution of the finite element discretization.
However, if the filter's radius is excessively large, it results in significant gray regions.
To overcome this, the filtered values, denoted by $\tilde{\mu}_e, e=1,\dots,N$, are translated to a near-binary density field using a smooth approximation of the Heaviside step function.

In this study, we employ this strategy, utilizing the commonly employed density filter on pseudo densities $\mu_e, e=1,\dots,N$, viz.,
\begin{align}
	\tilde{\mu}_e \coloneqq \frac{\sum_{i=1}^N w_{ei}  v_i \mu_i}{\sum_{i=1}^N w_{ei} v_i},
	\label{eq:density_filter}
\end{align}
where filter weights $w_{ei}$ are computed based on the distances between the finite element centroids $\boldsymbol{x}_e, e=1,\dots,N$, as
\begin{equation}
	w_{ei}        \coloneqq \max\{R -\Vert \boldsymbol{x}_e-\boldsymbol{x}_i \Vert, 0\},
	\label{eq:density_weights}
\end{equation}
with $R$ representing the filter radius.
Subsequently, the intermediate values $\tilde{\mu}_e, e=1,\dots,N$, undergo projection using a smooth approximation of the Heaviside step function to give material densities $\rho_e, e=1,\dots,N$, expressed as
\begin{align}
	\rho_e \coloneqq \frac{\tanh(\beta \eta) + \tanh(\beta (\tilde{\mu}_e - \eta))}{\tanh(\beta \eta) + \tanh(\beta (1 - \eta))}.
	\label{eq:heaviside}
\end{align}
The parameter $\beta > 0$ governs the sharpness of the projection, while $\eta \in [0, 1]$ acts as the threshold parameter.
Notably, when $\eta$ assumes values of $0$ and $1$, the projection corresponds to the Heaviside step filter~\citep{AchievingMinimGuest2004} and modified Heaviside filter~\citep{MorphologyBaseSigmun2007}, respectively.
In the limit where $\beta \rightarrow \infty$, these threshold values enable control over the length scale for the solid and the void phases, respectively.

\subsubsection{Surface grading of Young's modulus}

To incorporate surface grading of Young's modulus, we apply the nonlocal material model outlined in Section~\ref{sec:material_model} to the density field $\rho$.
For the approximation of the integral in Eq.~\eqref{eq:nonlocal_youngs_modulus}, we utilize mid-point integration, specifically one-point Gaussian quadrature.
Assuming the support of the kernel function $\psi(\,\cdot\,; \boldsymbol{x}_e)$ can be contained within the design domain for atleast one $\boldsymbol{x}_e, e=1,\dots,N$, the resulting integrated density $\bar{\rho}_e \in [0, 1]$ approximating $\phi$ is given by
\begin{align}
	\phi(\boldsymbol{x}_e) \approx \bar{\rho}_e \coloneqq \frac{\sum_{i=1}^N  \psi(\boldsymbol{x}_i; \boldsymbol{x}_e) v_i \rho_i}{\max_{j=1}^N{ \sum_{i=1}^N \psi(\boldsymbol{x}_i; \boldsymbol{x}_j) v_i}}.
	\label{eq:nonlocal_filter}
\end{align}
Here, the $\max$ operator in the denominator ensures that the nonlocal model behaves in a physically consistent manner at the edges of the domain.

Since $\bar{\rho}_e$ can be greater than zero even for void elements ($\rho_e=0$), we introduce a scalar field represented by $\alpha_e \in [0, 1], e=1,\dots,N$, where
\begin{equation}
	\alpha_e \coloneqq {\rho}_e^p \left[\zeta  + (1-\zeta) \bar{\rho}_e \right].
	\label{eq:elemental_youngs_modulus}
\end{equation}
The parameter $\alpha_e$ act as a scaling coefficient, which adjusts the reference Young's modulus $E_0$ to obtain the elemental Young's modulus.
The inclusion of the exponent $p > 1$ aims to penalize intermediate densities $\rho_e$, and encourages the emergence of a near-binary material distribution.
In conventional compliance minimization problems, a common choice is $p=3$ as used in the literature~\citep{TopologyOptimiBendso2004}.
However, due to the multiplication of the term $\rho_e^p$ by a linear function of the filtered density $\bar{\rho}_e$ in the above expression, lower values of $p$ depending on $\zeta$ are more appropriate to prevent excessive penalization that potentially slows down convergence in the optimization procedure.

The elemental Young's modulus corresponding to finite element $e$ is then expressed as
\begin{equation}
	E(\alpha_e) \coloneqq \left(\kappa + (1-\kappa) \alpha_e \right) E_0,
	\label{eq:elemental_youngs_modulus}
\end{equation}
where a minimal stiffness coefficient of $\kappa = 10^{-9}$ is incorporated for void elements ($\rho_e=\alpha_e=0$) to circumvent singularity issues in the global stiffness matrix during the numerical implementation.


\subsubsection{Optimization problem}

To exemplify the material model within the context of a structural design problem, we examine the standard compliance minimization problem with a structural volume constraint.
The optimization problem is formulated as follows:
\begin{equation}
	\begin{aligned}
		 & \min_{\boldsymbol{\mu}} \ c(\boldsymbol{\mu}) =  \mathbf{f}^\top \mathbf{u}(\boldsymbol{\alpha}(\boldsymbol{\mu})), \\
		 & \text{subject to }
		\begin{cases}
			g(\boldsymbol{\mu}) = \sum_{e=1}^N  v_e \rho_e(\boldsymbol{\mu}) \leq V, \\
			0 \leq \mu_e \leq 1, \quad e = 1,\dots,N.
		\end{cases}
	\end{aligned}
	\label{eq:optimization_problem}
\end{equation}
Here, $c : \mathbb{R}^N \rightarrow \mathbb{R}$ represents the global compliance, encompassing the overall stiffness of the structure.
The vectors $\mathbf{u} \in \mathbb{R}^M$ and $\mathbf{f} \in \mathbb{R}^M$ denote the nodal displacements and forces, respectively.
The structure's volume is defined by the function $g : \mathbb{R}^N \rightarrow \mathbb{R}_{\geq 0}$, and $V > 0$ is the prescribed volume.
For a given design $\boldsymbol{\mu}$ and external load $\mathbf{f}$, $\mathbf{u} = \mathbf{u}(\boldsymbol{\alpha}(\boldsymbol{\mu}))$ is the solution of the state problem:
\begin{equation}
	\mathbf{K}(\boldsymbol{\alpha}) \mathbf{u} = \mathbf{f},
	\text{ where }
	\mathbf{K}(\boldsymbol{\alpha}) = \sum_{e=1}^N E(\alpha_e) \mathbf{K}_e.
	\label{eq:state_equation}
\end{equation}
Here, $\mathbf{K}_e$ denotes the elemental stiffness matrix associated with element $e$, assuming a unit Young's modulus.


For solving the topology optimization problem (Eq.~\eqref{eq:optimization_problem}),
a solver based on the Method of Moving Asymptotes (MMA)~\citep{TheMethodOfMSvanbe1987} is used, implemented in the \texttt{ParOpt}~\citep{AScalableFramChin2019} library.
The gradient of the objective function as required by the optimizer is computed using the adjoint method~\citep{TopologyOptimiBendso2004}.

The sensitivity of the compliance $c$ with respect to the design parameter $\mu_e$ is expressed as
\begin{equation}
	\pdv{c}{\mu_e}
	=
	\sum_{i=1}^N
	\sum_{j=1}^N
	\pdv{c}{\alpha_i}
	\pdv{\alpha_i}{\rho_j}
	\pdv{\rho_j}{\tilde{\mu}_j}
	\pdv{\tilde{\mu}_j}{\mu_e},
\end{equation}
where the derivative $\partial c/\partial \alpha_e$ using the adjoint method is formulated as
\begin{equation}
	\pdv{c}{\alpha_e}
	=
	- E^\prime(\alpha_e) \mathbf{u}_e^\top \mathbf{K}_e \mathbf{u}_e.
\end{equation}
The other partial derivatives for $i,j=1,\dots,N$ are given as
\begin{align}
	\pdv{\alpha_i}{\rho_j}
	 & =
	p \rho_i^{p-1} \left[\zeta + (1-\zeta) \bar{\rho}_i \right] \delta_{ij}
	+ \rho_i^p (1-\zeta) \pdv{\bar{\rho}_i}{\rho_j},
	\\
	\pdv{\bar{\rho}_i}{\rho_j}
	 & =
	\frac{\psi(\boldsymbol{x}_j; \boldsymbol{x}_i) v_j}{\max_{l=1}^N{ \sum_{k=1}^N \psi(\boldsymbol{x}_k; \boldsymbol{x}_l) v_k}},
	\\
	\pdv{\rho_i}{\tilde{\mu}_i}
	 & =
	\frac{\beta [1-\tanh^2(\beta (\tilde{\mu}_i - \eta))]}{\tanh(\beta \eta) + \tanh(\beta (1 - \eta))},
	\\
	\pdv{\tilde{\mu}_i}{\mu_j}
	 & =
	\frac{w_{ij} v_j}{\sum_{k=1}^N w_{ik} v_k},
\end{align}
where $\delta_{ij}$ is Kronecker delta.
Similarly, using the chain rule, the derivative $\partial g/\partial\mu_e$ is given by
\begin{equation}
	\pdv{g}{\mu_e}
	=
	\sum_{i=1}^N \sum_{j=1}^N v_i \pdv{\rho_i}{\tilde{\mu}_j} \pdv{\tilde{\mu}_j}{\mu_e}.
\end{equation}

\section{Numerical examples}
\label{sec:results}

\begin{figure*}[htbp]
	\begin{center}
		\begin{tikzpicture}[scale=0.37]
			\scaling{0.37};
			\draw[thick, fill=black!10] (0, 0) rectangle (15, 10);
			\point{a}{0}{0};
			\point{b}{0}{10};
			\point{c}{0}{5};
			\point{a1}{-0.2}{0};
			\point{b1}{-0.2}{5};
			\foreach \y in {0,1,...,10} {
					\draw[thick,] (0, \y) -- (-0.5, \y-0.5);
				}
			\draw[fill=black!50] (0, 0) rectangle (0.5, 10);
			\draw[fill=black!50] (14.5, 4.5) rectangle (15, 5.5);

			\draw[ultra thick] (15.2, 4.5) -- (15.2, 5.5);
			\draw[-latex, ultra thick] (15.2, 5) --  (15.6, 5.) -- (15.6, 3) node[anchor=north] {$f$};
			\point{d}{15}{0};
			\dimensioning{2}{a}{b}{5.5}[\SI{10}{\milli\meter}];
			\dimensioning{1}{a}{d}{15}[\SI{15}{\milli\meter}];

			\point{a}{0}{0};
			\point{b}{0.5}{0};
			\dimensioning{1}{a}{b}{10.5}[\SI{0.5}{\milli\meter}];
			\point{a}{14.5}{4.5};
			\point{b}{15}{4.5};
			\dimensioning{1}{a}{b}{10.5}[\SI{0.5}{\milli\meter}];
			\point{a}{14.5}{4.5};
			\point{b}{14.5}{5.5};
			\dimensioning{2}{a}{b}{14}[\SI{1}{\milli\meter}];

			\begin{scope}[shift={(20, 0)}]
				\point{a}{0}{0};
				\point{a1}{0}{-0.1};
				\point{a2}{1}{-0.1};
				\point{b}{20}{10};
				\point{f}{20}{9.5};
				\point{c}{20}{0};
				\point{g}{20}{0.5};
				\point{d}{0}{5};
				\point{e}{1}{0.5};
				\draw[fill=black!10, thick] (a) rectangle (b);
				\draw[fill=black!50] (a) rectangle (e);
				\draw[fill=black!50] (19, 0) rectangle (20, 0.5);
				\point{b1}{19}{-0.1};
				\point{b2}{20}{-0.1};
				\foreach \i in {0,1,...,20} {
						\draw[-latex, thick] (\i, 12) -- (\i, 10.2);
					}
				\draw[fill=black!50] (0, 9.5) rectangle (20, 10);
				\dimensioning{1}{a}{c}{15}[$\SI{20}{\milli\meter}$];
				\dimensioning{2}{f}{b}{-1}[$\SI{0.5}{\milli\meter}$];
				\dimensioning{2}{c}{g}{-1}[$\SI{0.5}{\milli\meter}$];
				\dimensioning{1}{b1}{b2}{13}[$\SI{1}{\milli\meter}$];
				\dimensioning{1}{a1}{a2}{13}[$\SI{1}{\milli\meter}$];
				\foreach \i in {0,0.5,...,1} {
						\draw[thick] (\i, -0.) -- (\i-0.5, -0.5);
					}
				\foreach \i in {19,19.5,...,20} {
						\draw[thick] (\i, -0.) -- (\i-0.5, -0.5);
					}
				\node[right] at (16., 18) {solid};
				\draw[fill=black!50] (15, 17.5) rectangle (16, 18.5);
				\node[right] at (11., 18) {design};
				\draw[fill=black!10] (10, 17.5) rectangle (11, 18.5);
				\node at (10, 13) {$f$};
			\end{scope}
		\end{tikzpicture}
	\end{center}
	\caption{
		Geometrical setup of 2D cantilever beam (left) and bridge (right) structures.
	}
	\label{fig:design_domain_schematic}
\end{figure*}
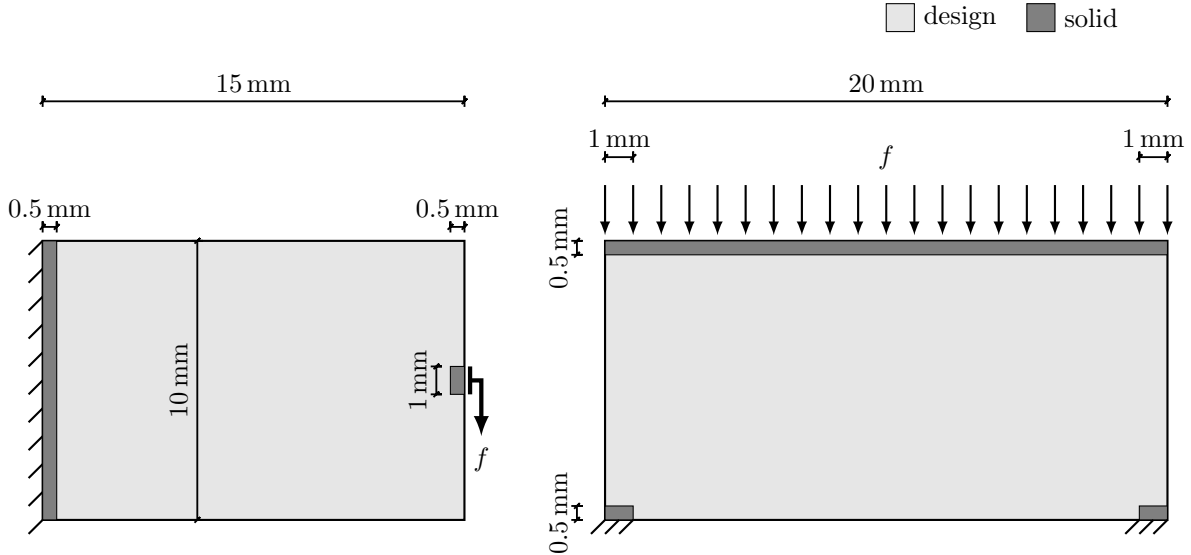

\newlength{\imagewidth}
\settowidth{\imagewidth}{\includegraphics{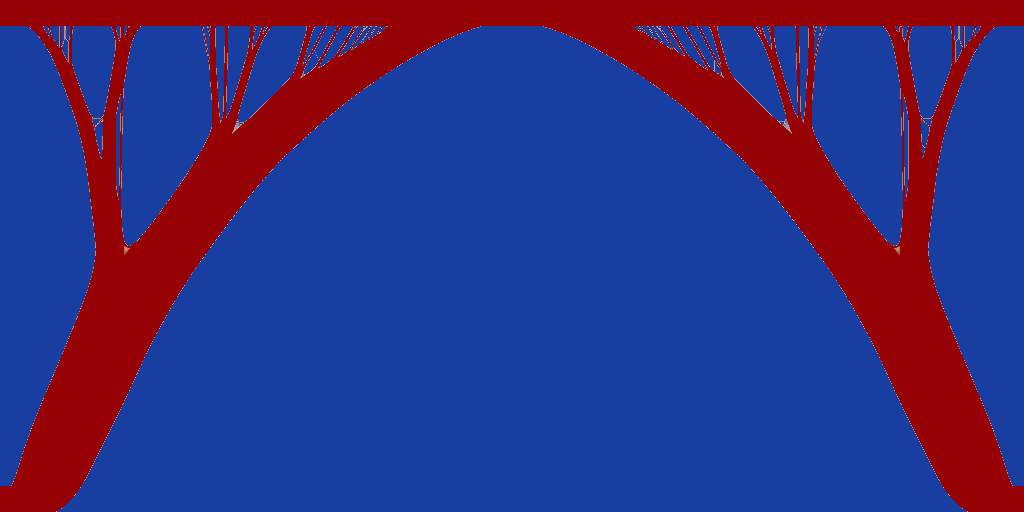}}
\begin{figure*}[htbp]
	\flushright{
		\includegraphics[width=0.25\textwidth]{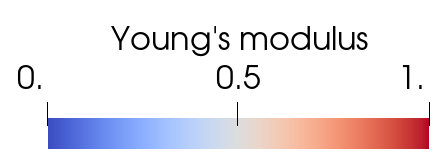}
	}
	\begin{subfigure}{\textwidth}
		\centering
		\begin{tabular}[c]{cc|cccc}
			 &
			 & \multicolumn{4}{c}{Cantilever beam -- domain size}
			\\[0.2cm]
			 &
			 & $15\times10$ mm$^2$
			 & $30\times20$ mm$^2$
			 & $45\times30$ mm$^2$
			 & $60\times40$ mm$^2$
			\\
			\hline
			\\
			\multirow{4}{*}{\rotatebox[origin=c]{90}{Optimized with}}
			 & \rotatebox[origin=l]{90}{ \begin{tabular}{c} $\delta=0$ \end{tabular} }
			 & \includegraphics[width=0.2\textwidth]{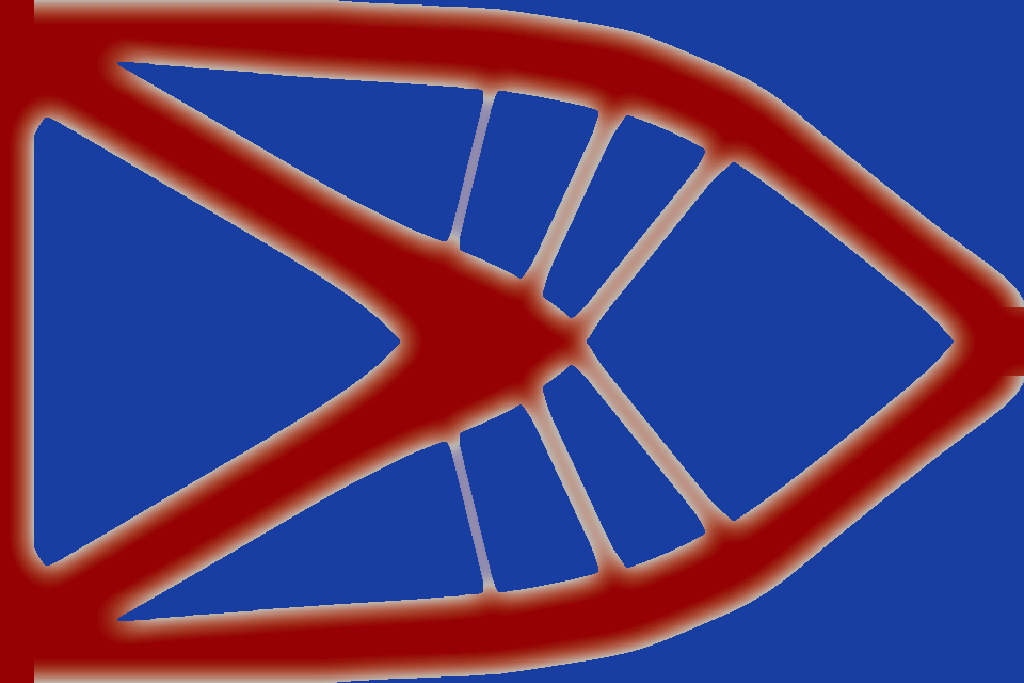}
			 & \includegraphics[width=0.2\textwidth]{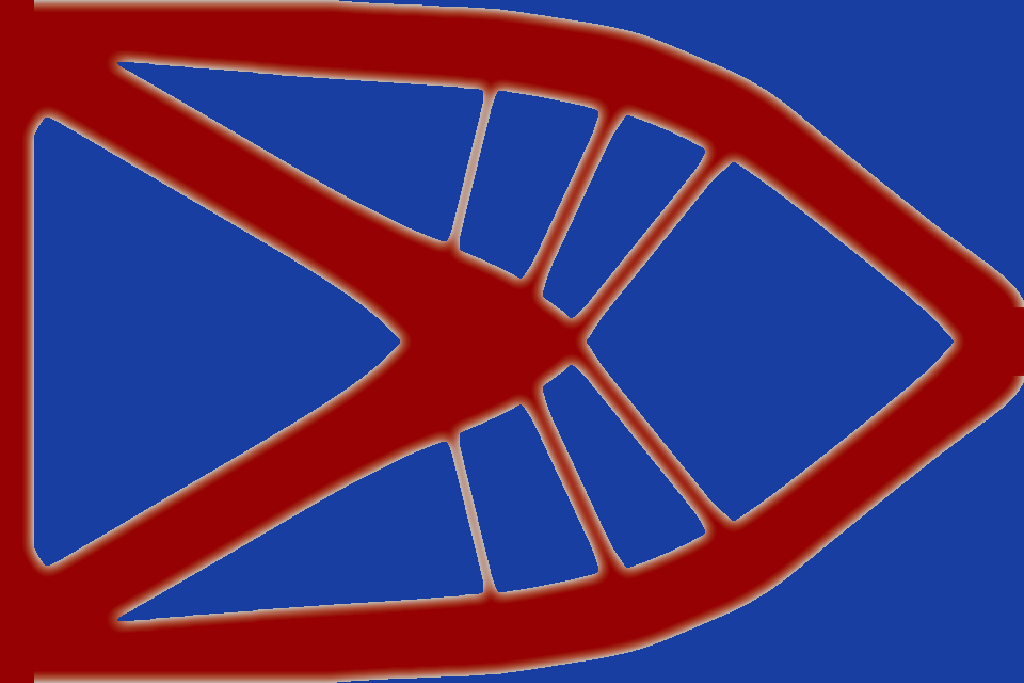}
			 & \includegraphics[width=0.2\textwidth]{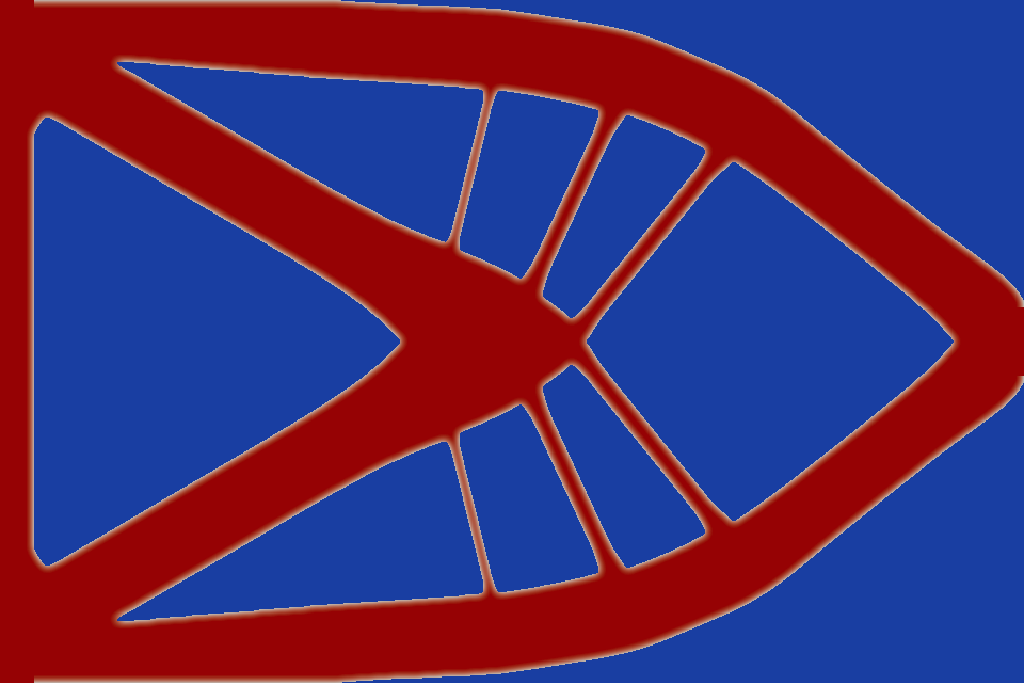}
			 & \includegraphics[width=0.2\textwidth]{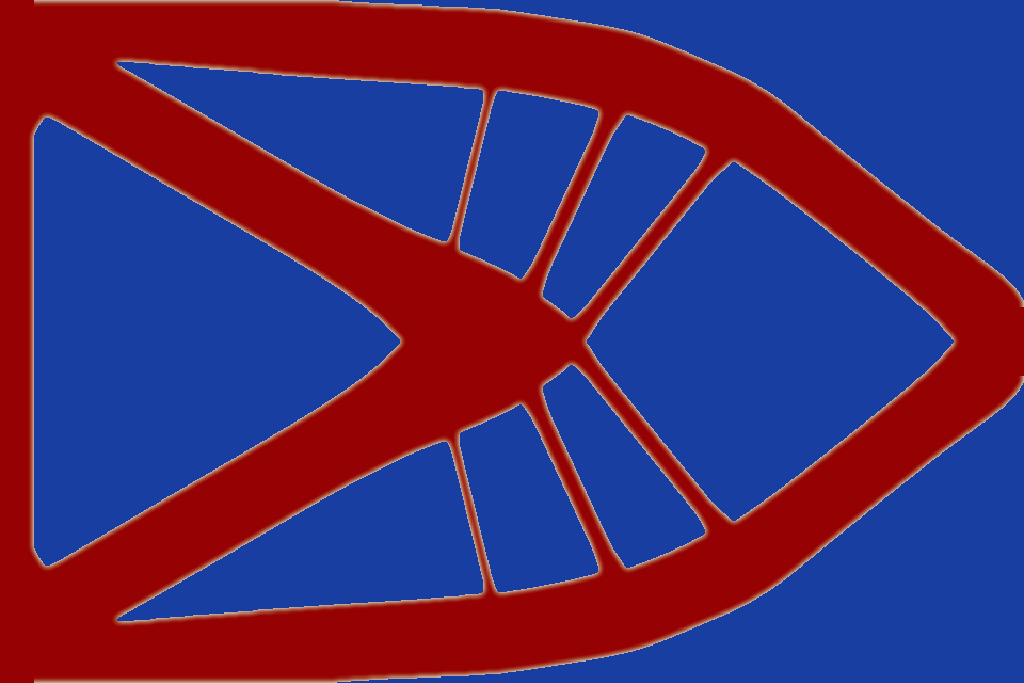}
			\\
			 &
			 & $c=\SI{39.49}{\newton\milli\meter}$
			 & $c=\SI{146.37}{\newton\milli\meter}$
			 & $c=\SI{321.73}{\newton\milli\meter}$
			 & $c=\SI{565.59}{\newton\milli\meter}$
			\\[0.2cm]
			 & \rotatebox[origin=l]{90}{ \begin{tabular}{c} $\delta=\SI{0.4}{\milli\meter}$ \end{tabular}}
			 & \includegraphics[width=0.2\textwidth]{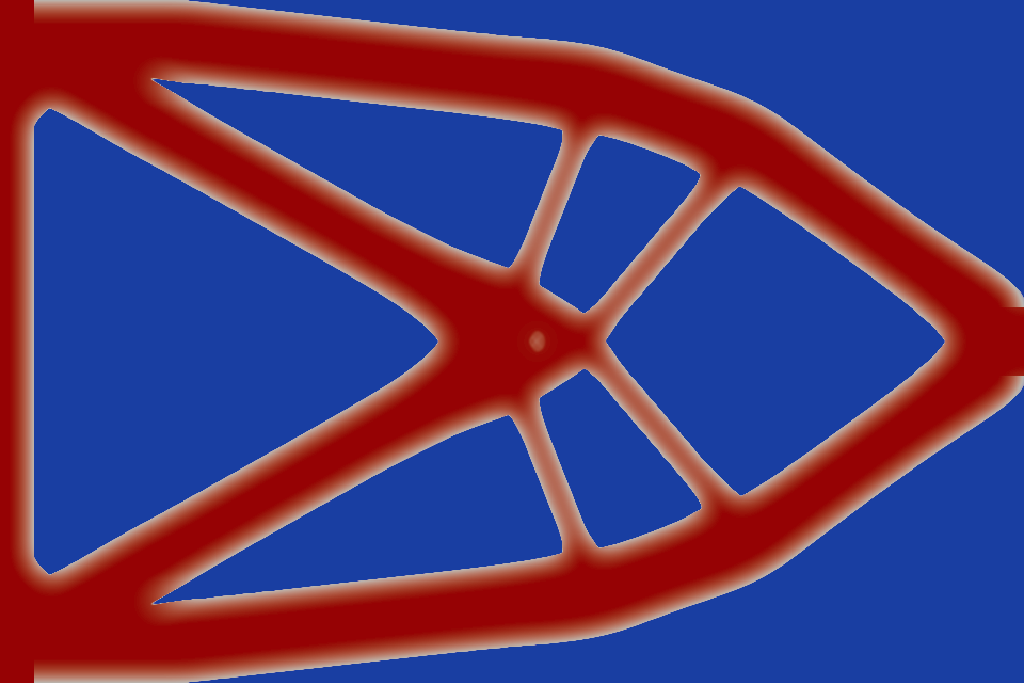}
			 & \includegraphics[width=0.2\textwidth]{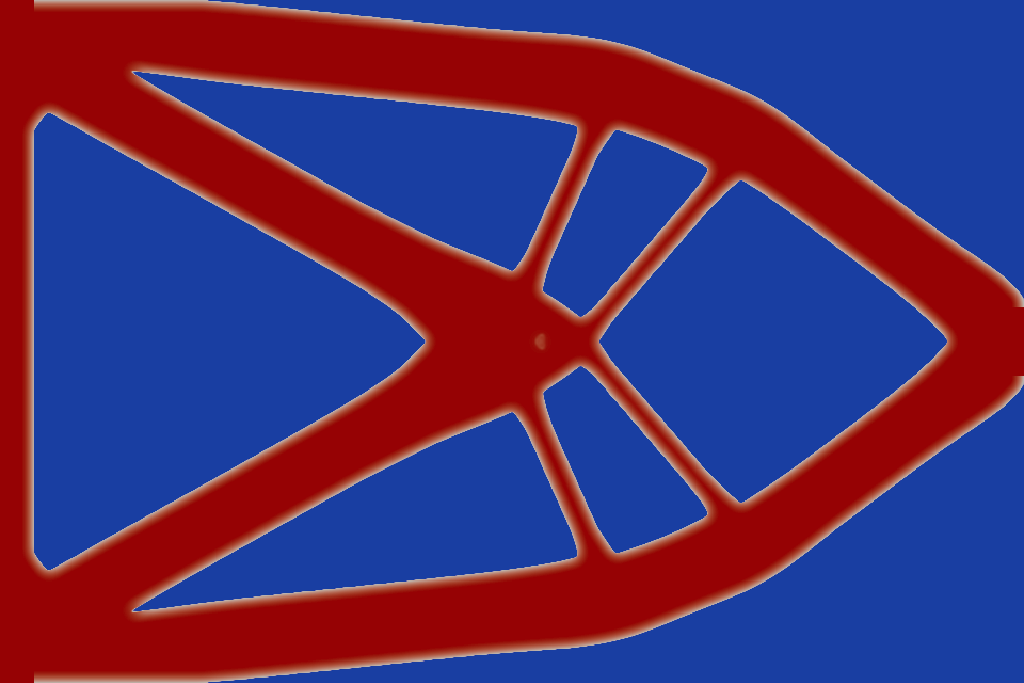}
			 & \includegraphics[width=0.2\textwidth]{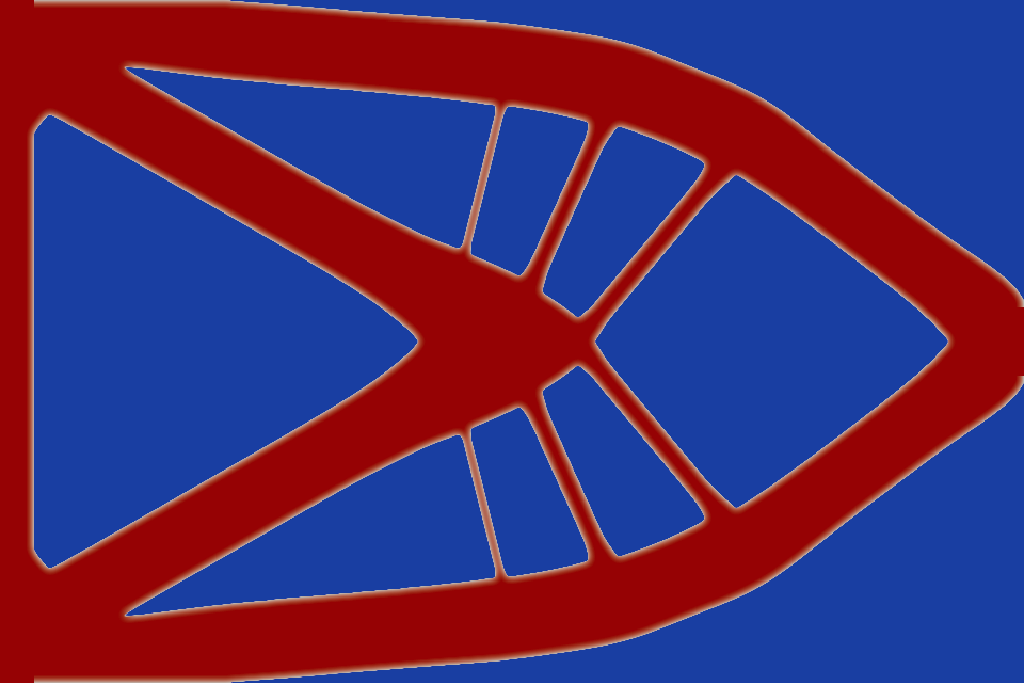}
			 & \includegraphics[width=0.2\textwidth]{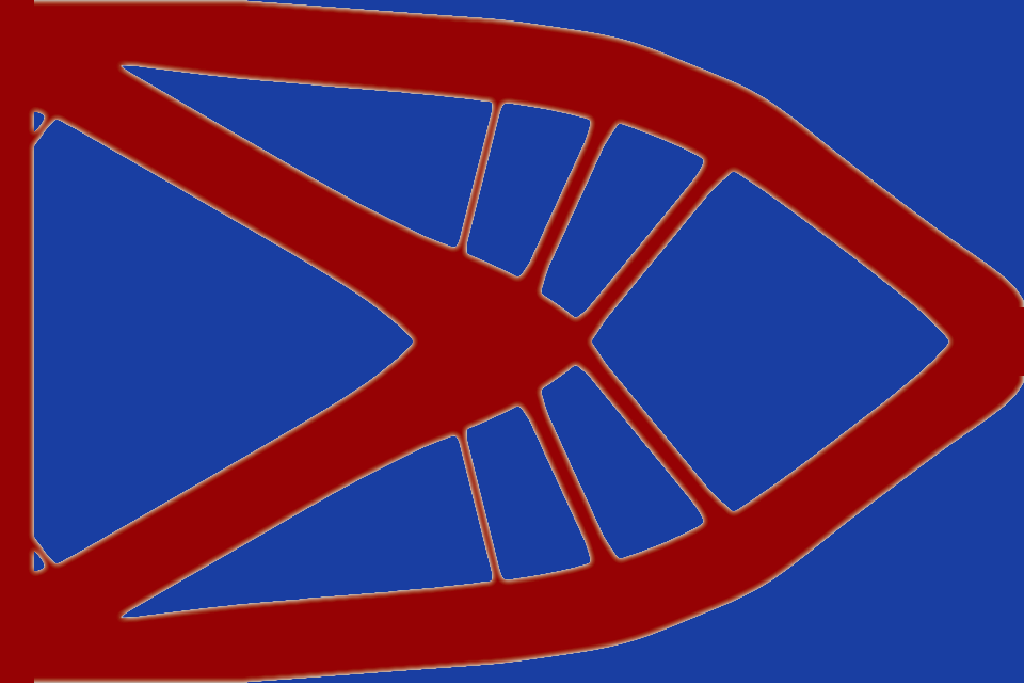}
			\\
			 &
			 & $c=\SI{38.28}{\newton\milli\meter}$
			 & $c=\SI{144.66}{\newton\milli\meter}$
			 & $c=\SI{320.04}{\newton\milli\meter}$
			 & $c=\SI{563.31}{\newton\milli\meter}$
			\\
		\end{tabular}
	\end{subfigure}
	\vskip 0.3in
	\begin{subfigure}{\textwidth}
		\centering
		\begin{tabular}[c]{cc|cccc}
			 &
			 & \multicolumn{4}{c}{Bridge -- domain size}
			\\[0.2cm]
			 &
			 & $20\times10$ mm$^2$
			 & $40\times20$ mm$^2$
			 & $60\times30$ mm$^2$
			 & $80\times40$ mm$^2$
			\\
			\hline
			\\
			\multirow{4}{*}{\rotatebox[origin=c]{90}{Optimized with \quad \quad }}
			 & \rotatebox[origin=c]{90}{\begin{tabular}{c} $\delta=0$ \end{tabular} }
			 & \raisebox{-.5\height}{\includegraphics[trim=0.5\imagewidth{} 0 0 0, clip, width=0.2\textwidth]{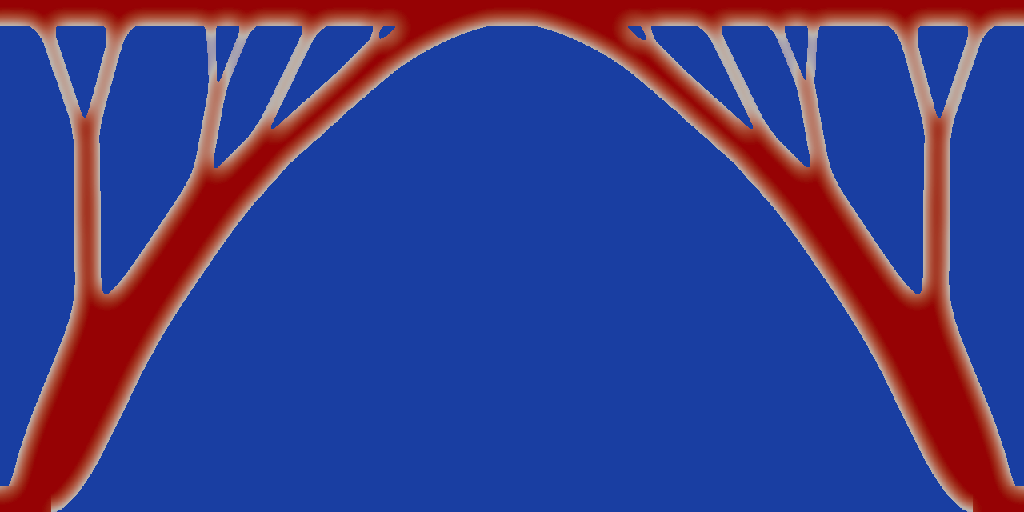}}
			 & \raisebox{-.5\height}{\includegraphics[trim=0.5\imagewidth{} 0 0 0, clip, width=0.2\textwidth]{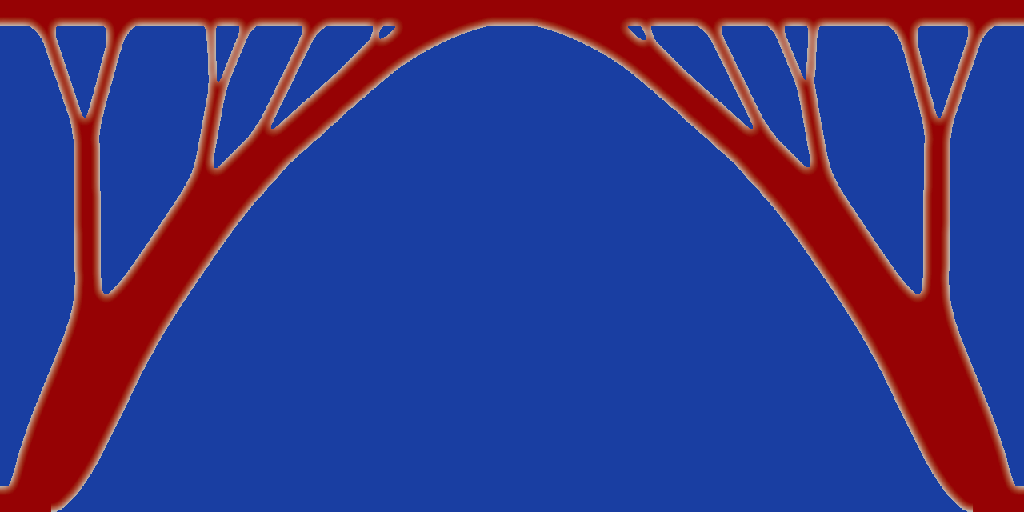}}
			 & \raisebox{-.5\height}{\includegraphics[trim=0.5\imagewidth{} 0 0 0, clip, width=0.2\textwidth]{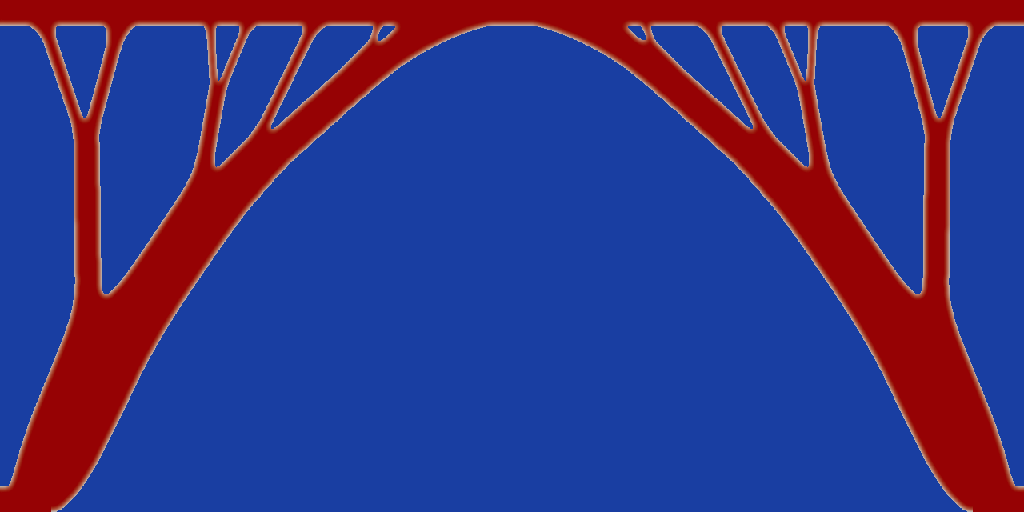}}
			 & \raisebox{-.5\height}{\includegraphics[trim=0.5\imagewidth{} 0 0 0, clip, width=0.2\textwidth]{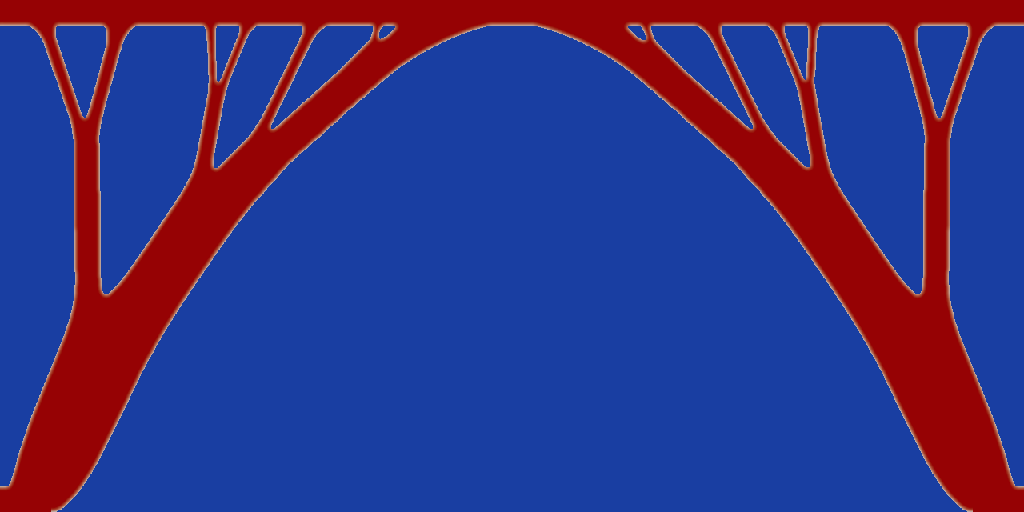}}
			\\
			 &
			 & $c=\SI{1985.19}{\newton\milli\meter}$
			 & $c=\SI{7177.46}{\newton\milli\meter}$
			 & $c=\SI{15655.2}{\newton\milli\meter}$
			 & $c=\SI{27420.1}{\newton\milli\meter}$
			\\[0.4cm]
			 & \rotatebox[origin=c]{90}{\begin{tabular}{c} $\delta=\SI{0.4}{\milli\meter}$ \end{tabular}}
			 & \raisebox{-.5\height}{\includegraphics[trim=0.5\imagewidth{} 0 0 0, clip, width=0.2\textwidth]{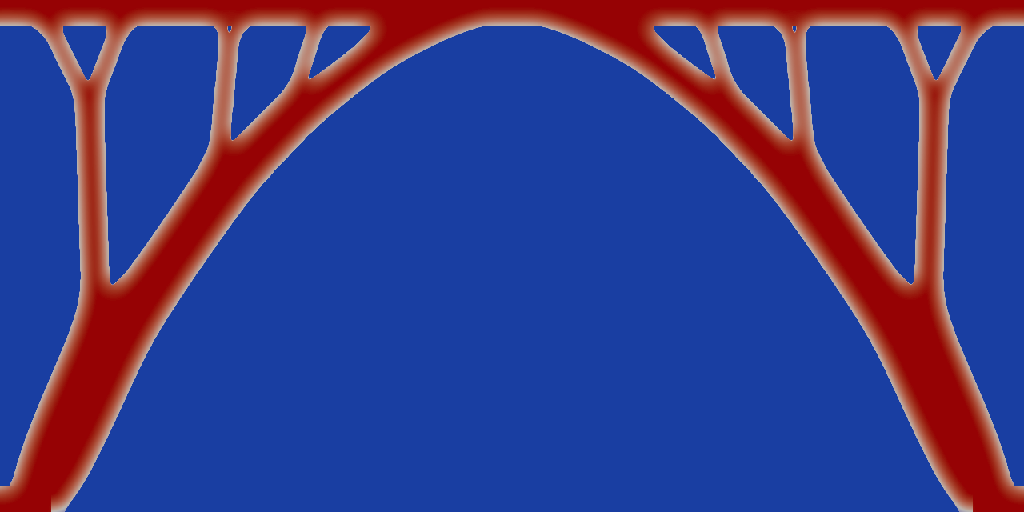}}
			 & \raisebox{-.5\height}{\includegraphics[trim=0.5\imagewidth{} 0 0 0, clip, width=0.2\textwidth]{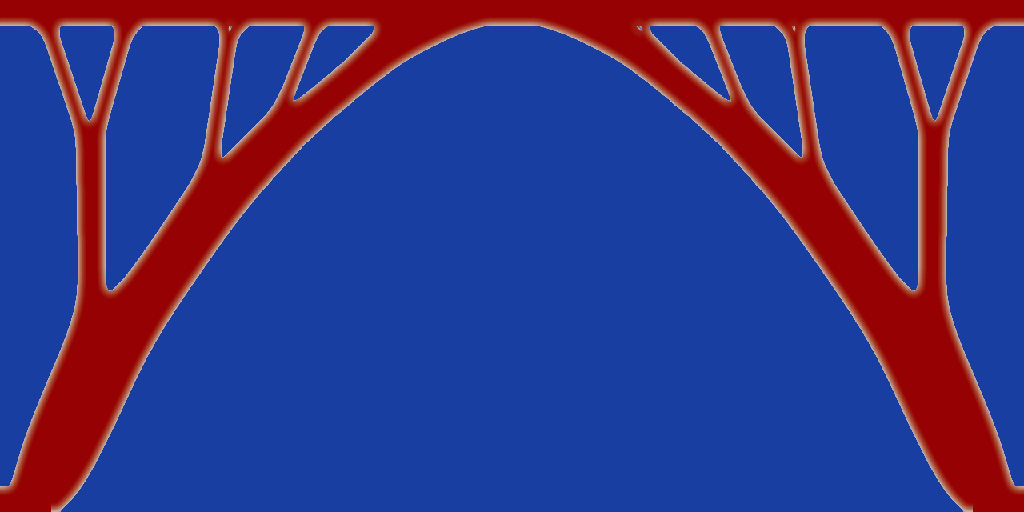}}
			 & \raisebox{-.5\height}{\includegraphics[trim=0.5\imagewidth{} 0 0 0, clip, width=0.2\textwidth]{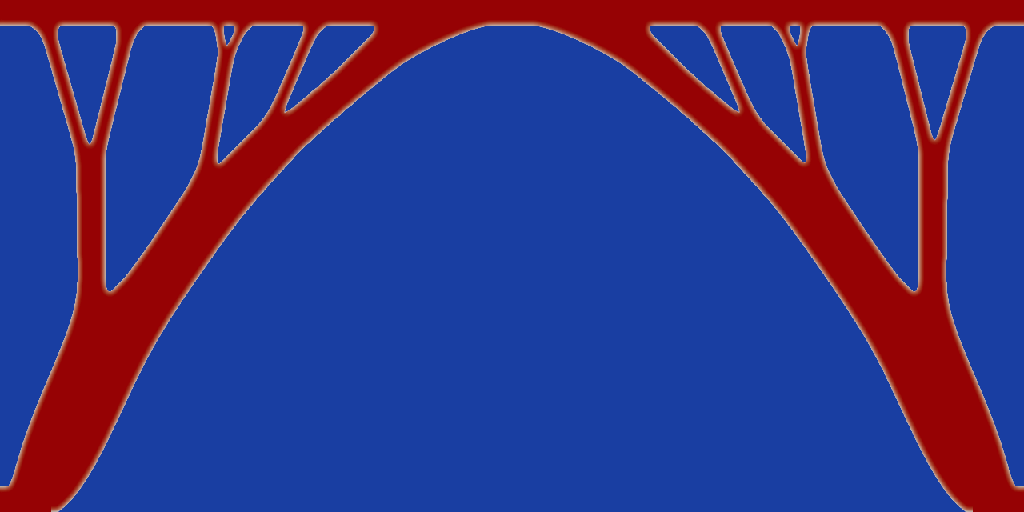}}
			 & \raisebox{-.5\height}{\includegraphics[trim=0.5\imagewidth{} 0 0 0, clip, width=0.2\textwidth]{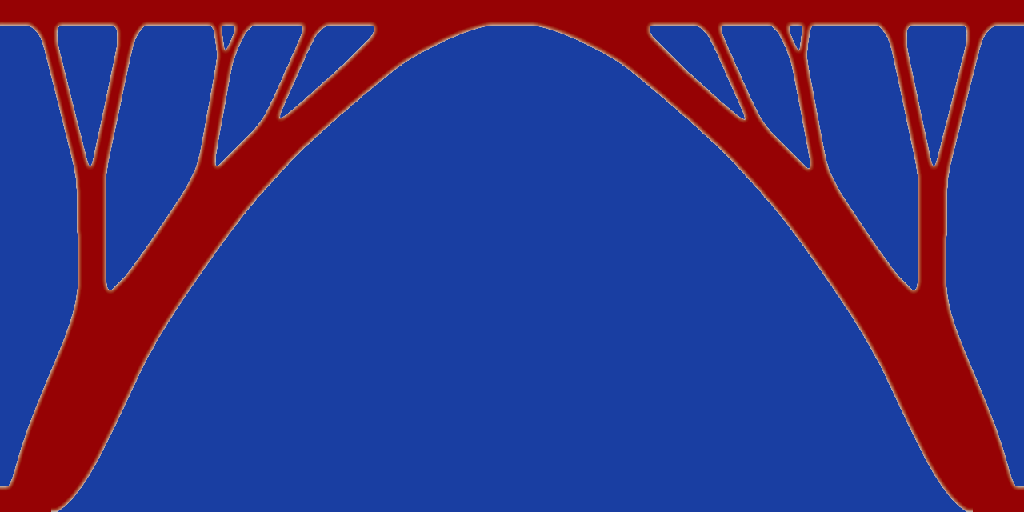}}
			\\
			 &
			 & $c=\SI{1942.21}{\newton\milli\meter}$
			 & $c=\SI{7152.23}{\newton\milli\meter}$
			 & $c=\SI{15645.8}{\newton\milli\meter}$
			 & $c=\SI{27405.0}{\newton\milli\meter}$
		\end{tabular}
	\end{subfigure}
	\vskip 0.2in
	\caption{
		Final designs for the 2D cantilever beam and bridge structures, optimized with $\delta=0$ and $\delta=\SI{0.4}{\milli\meter}$, and varying sizes of the design domain.
		For all designs, the Young's modulus and compliance values are computed in a postprocessing step where a volume preserving thresholding is performed and nonlocal material model is applied with $\delta=\SI{0.4}{\milli\meter}$.
		Owing to the bridge's symmetry across the central vertical axis, only the right half of the structure is depicted.
	}
	\label{fig:optimized_designs_scaling}
\end{figure*}

In the following section, we apply the presented optimization framework for the design of cantilever beam and bridge structures in two dimensions.
A plain-strain condition is assumed.
We set the reference Young's modulus to $E_0=\SI{1}{\mega\pascal}$, and the Poisson's ratio is assumed to be \num{0.3}.

The geometrical dimensions of the two structures are illustrated in Figure~\ref{fig:design_domain_schematic}.
The domain is discretized using bilinear quadrilateral elements with thickness of \SI{1}{\milli\meter}.
The cantilever beam is fixed at the left end, while a static traction of \SI{1}{\mega\pascal} distributed over \SI{1}{\milli\meter} length is applied at the right end.
For the bridge structure, a distributed load of \SI{1}{\mega\pascal} is applied on the top surface, while the bottom corners are held fixed.
The regions in dark gray indicate non-design regions, introduced to have full elastic stiffness at the support and loading points.
The volume fraction for the design domain is set to $\num{0.4}$ for the cantilever beam problem and $\num{0.2}$ for the bridge problem.

\subsection{Size effect}

\begin{figure*}[htbp]
	\centering
	\begin{subfigure}[b]{0.49\textwidth}
		\centering
		\includegraphics[width=\textwidth]{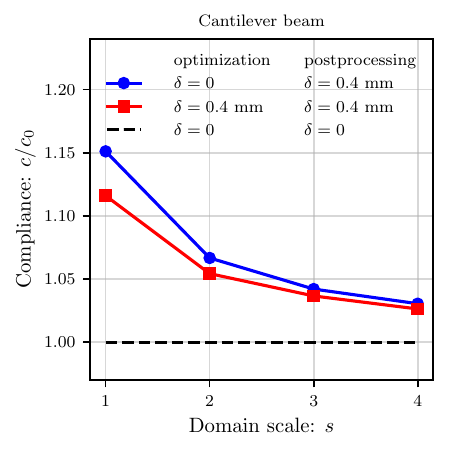}
	\end{subfigure}
	\begin{subfigure}[b]{0.49\textwidth}
		\centering
		\includegraphics[width=\textwidth]{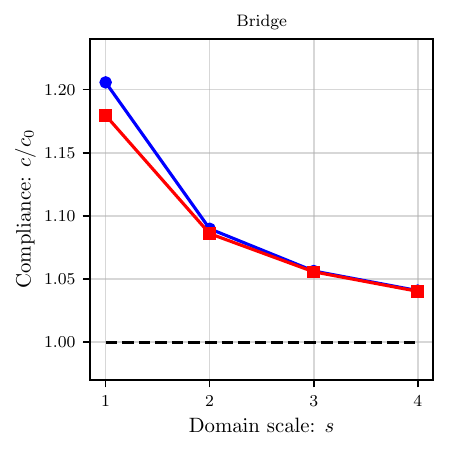}
	\end{subfigure}
	\caption{
		Variation in structural compliance with respect to the domain size for the cantilever beam and bridge structures, where the geometric dimensions are scaled by a factor of $s$, while keeping $\delta$ fixed.
		The compliance values ($c$) are normalized by the compliance ($c_0$) of the designs optimized and postprocessed with $\delta=0$.
		Designs optimized with surface grading effects ($\delta=\SI{0.4}{\milli\meter}$) exhibit lower compliance values in comparison to those generated without the integration of surface grading ($\delta=0$) in the optimization process.
	}
	\label{fig:scale_compliance}
\end{figure*}

First, we examine the influence of the overall design domain size on the optimized designs and their corresponding structural compliance.
Optimizations are conducted on distinct domain sizes, scaled by factors of $s=1, 2, 3$ and $4$.
These correspond to domain dimensions of $15\times10$ mm$^2, 30\times20$ mm$^2, 45\times30$ mm$^2,$ and $60\times40$ mm$^2$ for the cantilever beam.
For the bridge structure, the domain sizes are $20\times10$ mm$^2, 40\times20$ mm$^2, 60\times30$ mm$^2,$ and $80\times40$ mm$^2$, respectively.

Despite variations in domain size, the finite element discretization is consistently executed with the same number of elements.
Consequently, the mesh size $h$ is proportionally scaled by the factor $s$.
For the coarsest mesh in both the cantilever beam and bridge cases, $h=H=\SI{1.5625e-2}{\milli\meter}$.
The radius of the Heaviside projection filter is set to $R=10h$.
The continuation parameter $\beta$ is systematically increased at intervals of every \num{50} optimization iterations, taking values of $\beta=1, 2, 4, 8, 16, 32$ and  $64$.
After reaching the highest value of $\beta$, the optimization process is continued until convergence, with a maximum of \num{1000} iterations.
The threshold parameter $\eta$ is chosen to be \num{0.5}.

For optimization, two distinct values of the kernel parameter $\delta$ are considered. The first value, $\delta=0$, corresponds to the standard SIMP approach with the Heaviside projection filter. Throughout the optimization, the penalization constant $p$ remains fixed at $p=3$. The second value, $\delta=\SI{0.4}{\milli\meter}$, enables the emulation of surface grading effects. For this case, the penalization constant is selected as $p=2.5$. The value of $\zeta$ is set to $0$.


Figure~\ref{fig:optimized_designs_scaling} showcases the optimized designs.
The Young's modulus distribution and compliance values depicted are derived from a postprocessing stage involving volume-preserving thresholding followed by the employment of nonlocal material model with $\delta=\SI{0.4}{\milli\meter}$.
It is noteworthy that the designs incorporating surface grading exhibit significantly lower compliance compared to those generated without the explicit consideration of surface grading in the optimization procedure.
Moreover, as the domain size increases, the influence of surface grading diminishes, leading to nearly identical designs for both $\delta$ values.
Figure~\ref{fig:scale_compliance} shows the relationship between domain scale and compliance, with compliance values normalized by those obtained for designs optimized and postprocessed with $\delta=0$.
This normalization enables direct comparison of the compliance performance across different domain sizes.
The results clearly demonstrate that the normalized compliance values increase with decreasing domain size, indicating a heightened sensitivity to optimizing smaller-scale structures.
Additionally, the difference between the compliance values optimized with and without surface grading effects grows with decreasing domain size, further emphasizing the beneficial impact of incorporating material information into the optimization process.

\begin{figure*}[htbp]
	\centering
	\begin{tabular}[c]{c|cc|cc}
		\multirow{2}{*}{\begin{tabular}{c} \\ Mesh \\ size \\ ($h$) \end{tabular}}
		 & \multicolumn{2}{c}{Optimized with $\delta=0$}
		 & \multicolumn{2}{c}{Optimized with $\delta=\SI{0.4}{\milli\meter}$}
		\\[0.3cm]
		 & \begin{tabular}{c} Material \\ distribution \\ ($\rho$) \end{tabular}
		 & \begin{tabular}{c} Young's modulus \\ ($\alpha$) \end{tabular}
		 & \begin{tabular}{c} Material \\ distribution \\ ($\rho$) \end{tabular}
		 & \begin{tabular}{c} Young's modulus \\ ($\alpha$) \end{tabular}
		\\
		\hline
		\\
		$H$
		 & \raisebox{-.5\height}{\includegraphics[width=0.2\textwidth]{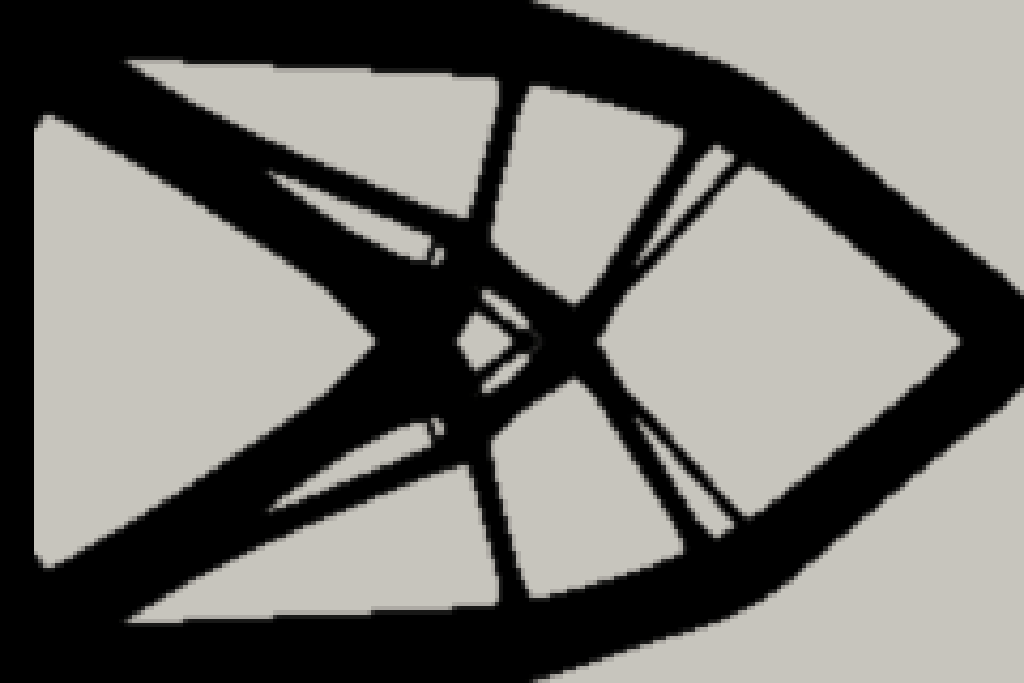}}
		 & \raisebox{-.5\height}{\includegraphics[width=0.2\textwidth]{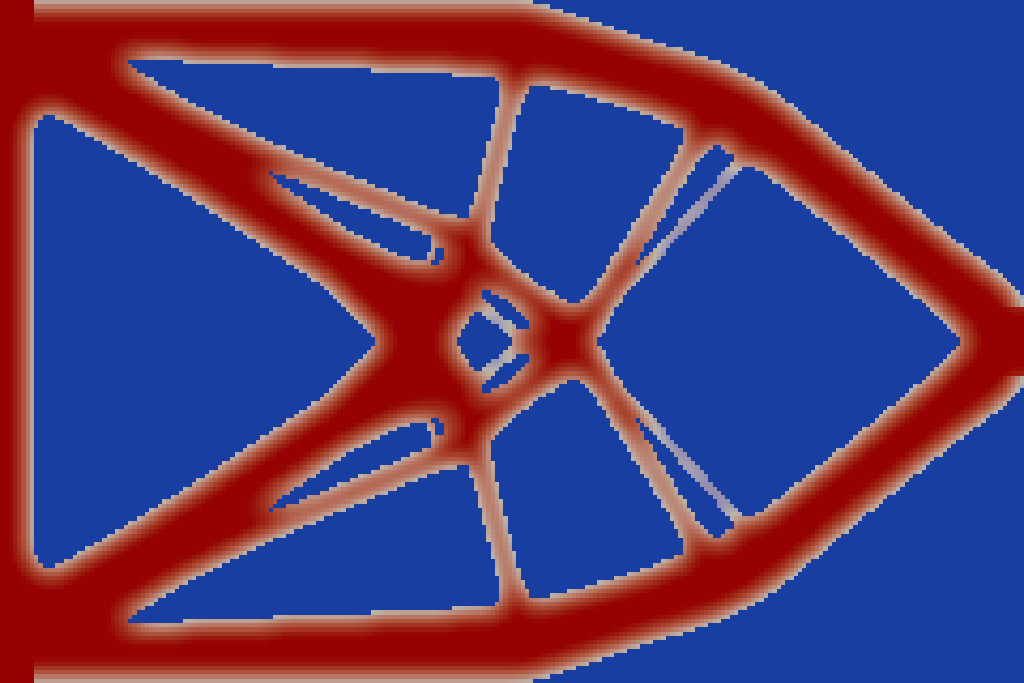}}
		 & \raisebox{-.5\height}{\includegraphics[width=0.2\textwidth]{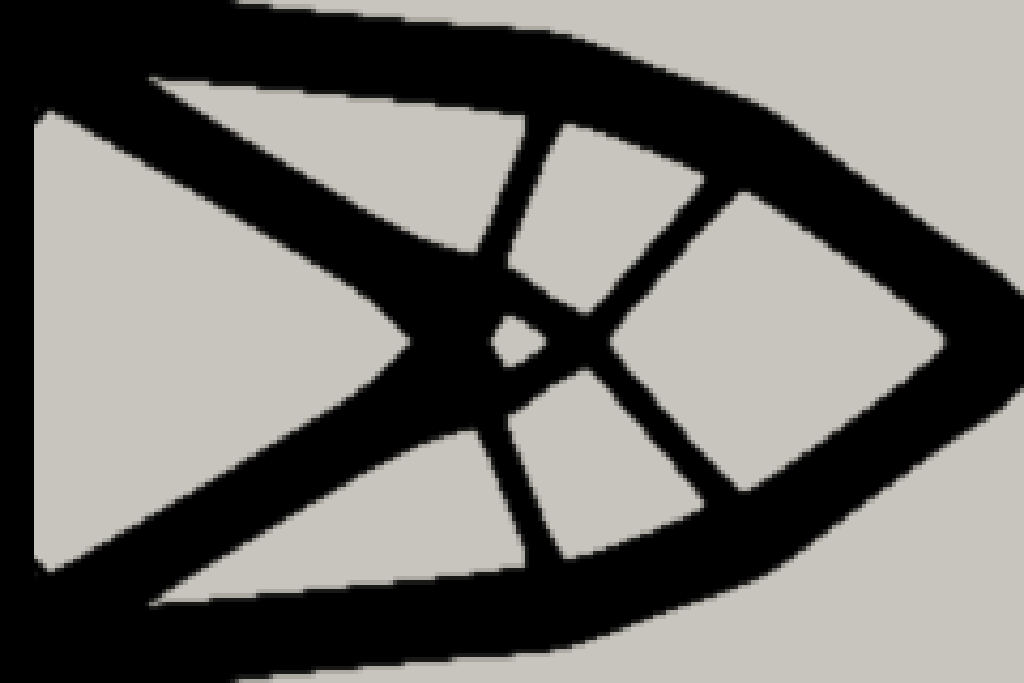}}
		 & \raisebox{-.5\height}{\includegraphics[width=0.2\textwidth]{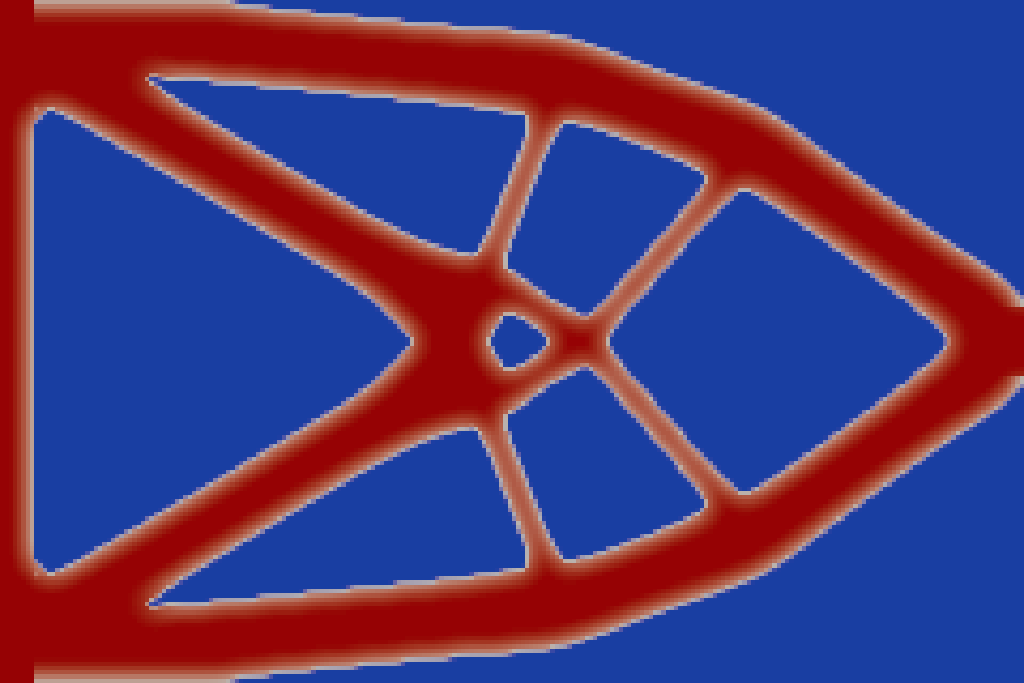}}
		\\
		 &
		 & $c=\SI{88.32}{\newton\milli\meter}$
		 &
		 & $c=\SI{83.63}{\newton\milli\meter}$
		\\[4mm]
		$\dfrac{H}{2}$
		 & \raisebox{-.5\height}{\includegraphics[width=0.2\textwidth]{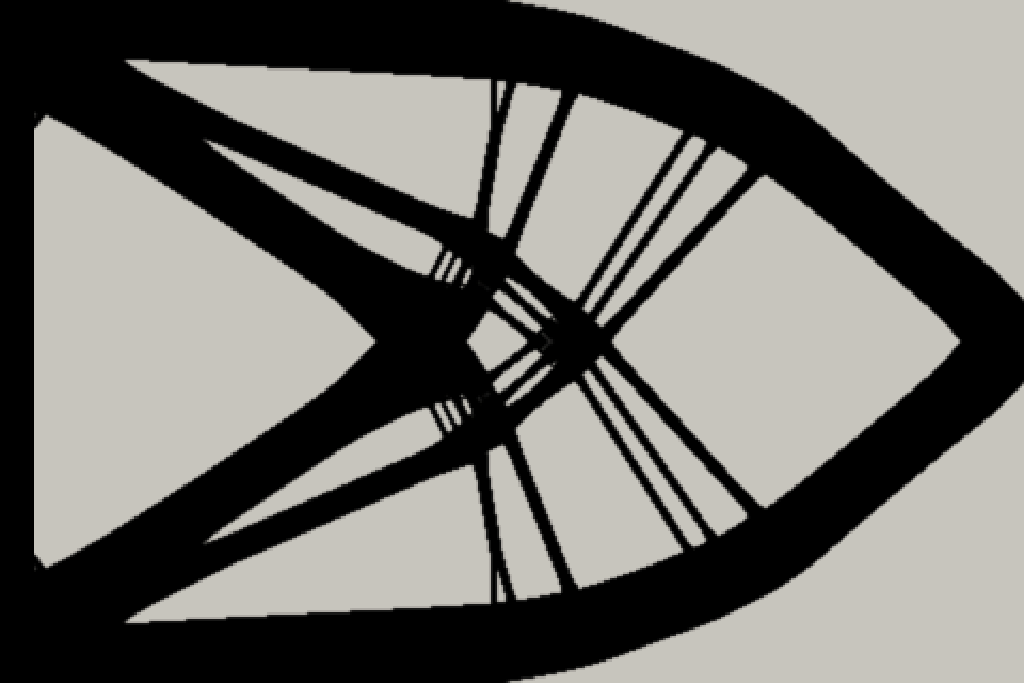}}
		 & \raisebox{-.5\height}{\includegraphics[width=0.2\textwidth]{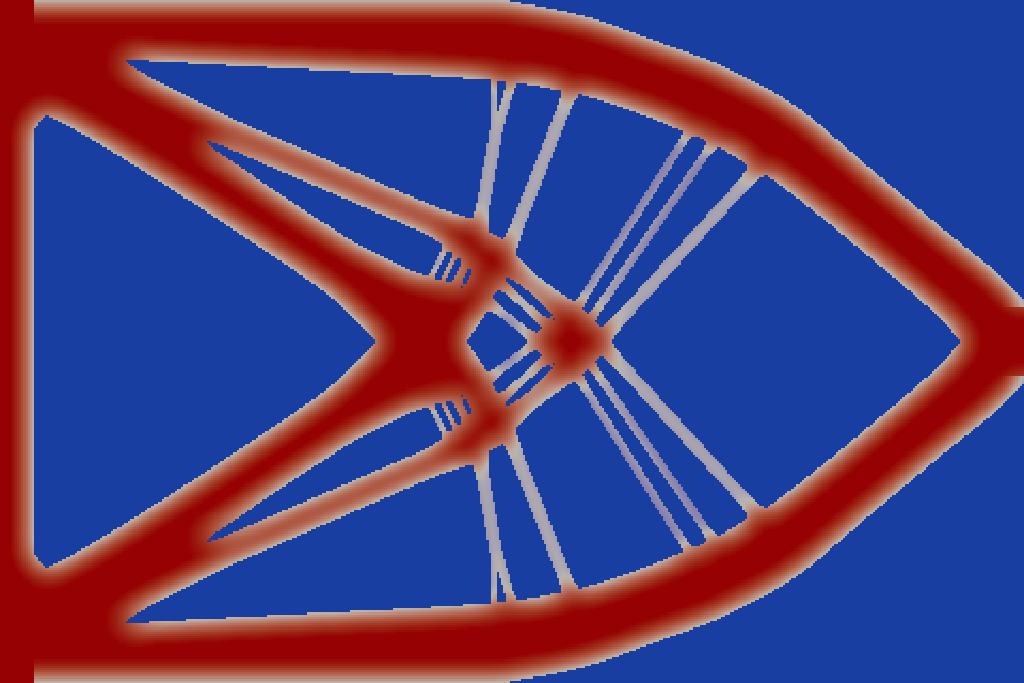}}
		 & \raisebox{-.5\height}{\includegraphics[width=0.2\textwidth]{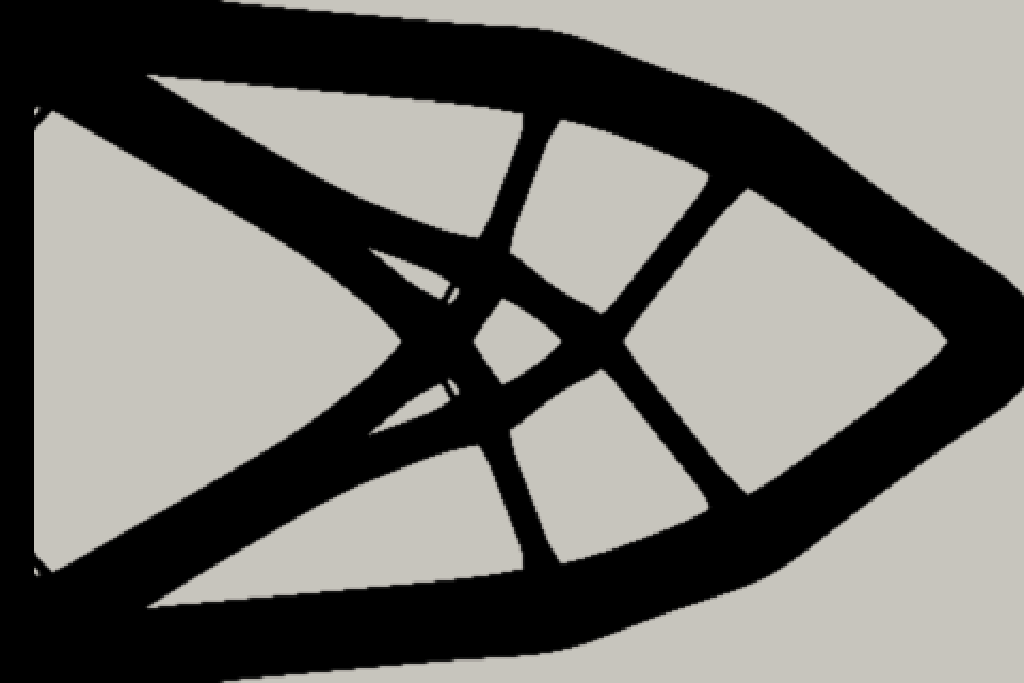}}
		 & \raisebox{-.5\height}{\includegraphics[width=0.2\textwidth]{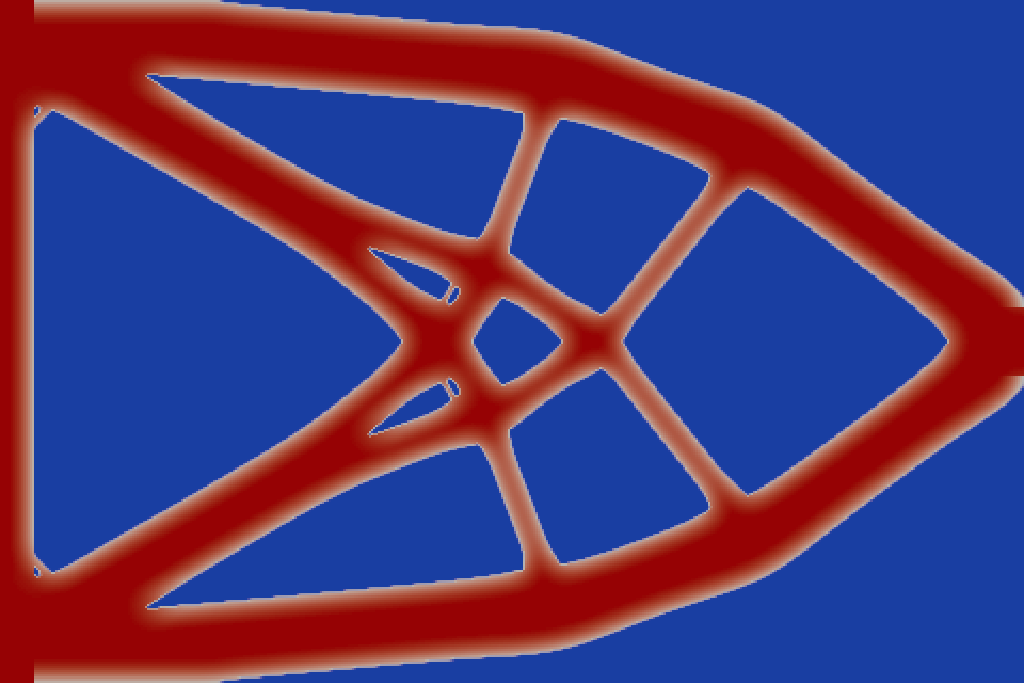}}
		\\
		 &
		 & $c=\SI{91.43}{\newton\milli\meter}$
		 &
		 & $c=\SI{83.20}{\newton\milli\meter}$
		\\[4mm]
		$\dfrac{H}{4}$
		 & \raisebox{-.5\height}{\includegraphics[width=0.2\textwidth]{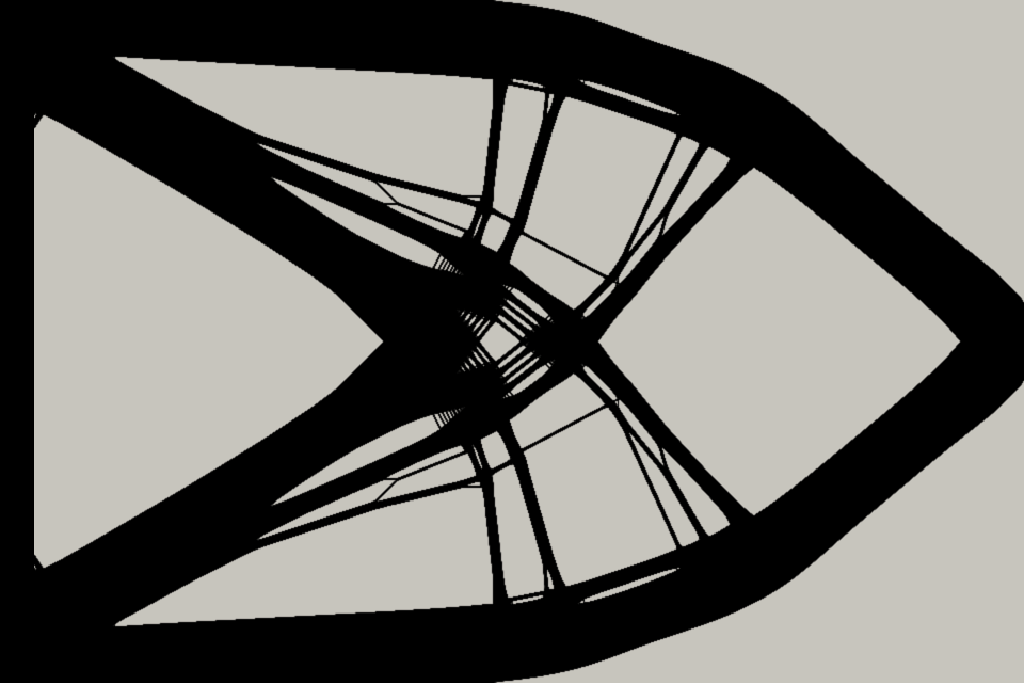}}
		 & \raisebox{-.5\height}{\includegraphics[width=0.2\textwidth]{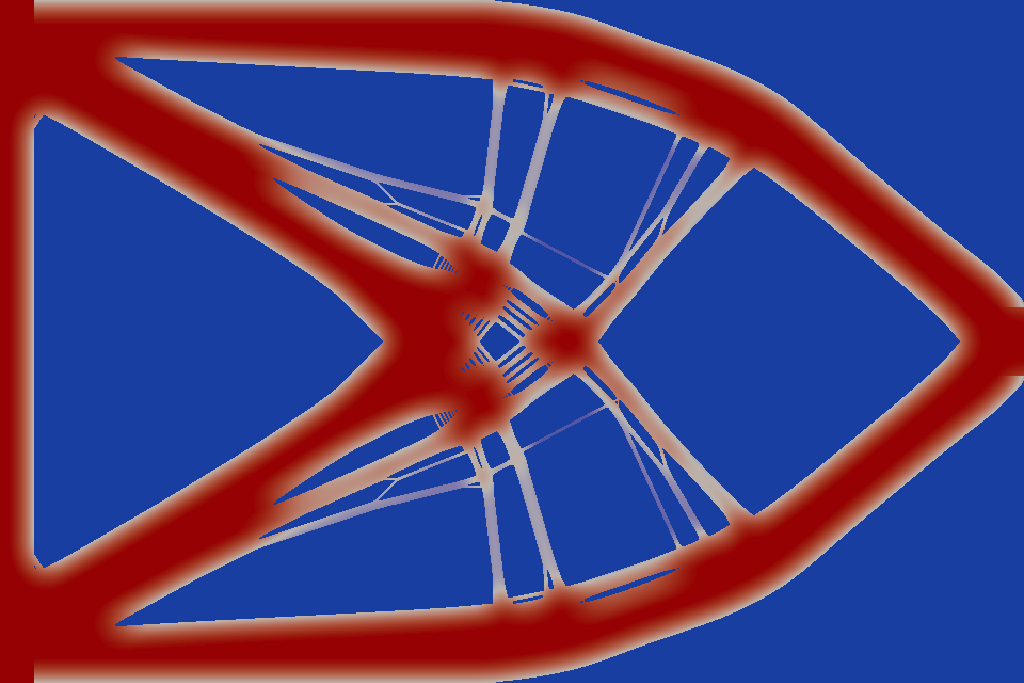}}
		 & \raisebox{-.5\height}{\includegraphics[width=0.2\textwidth]{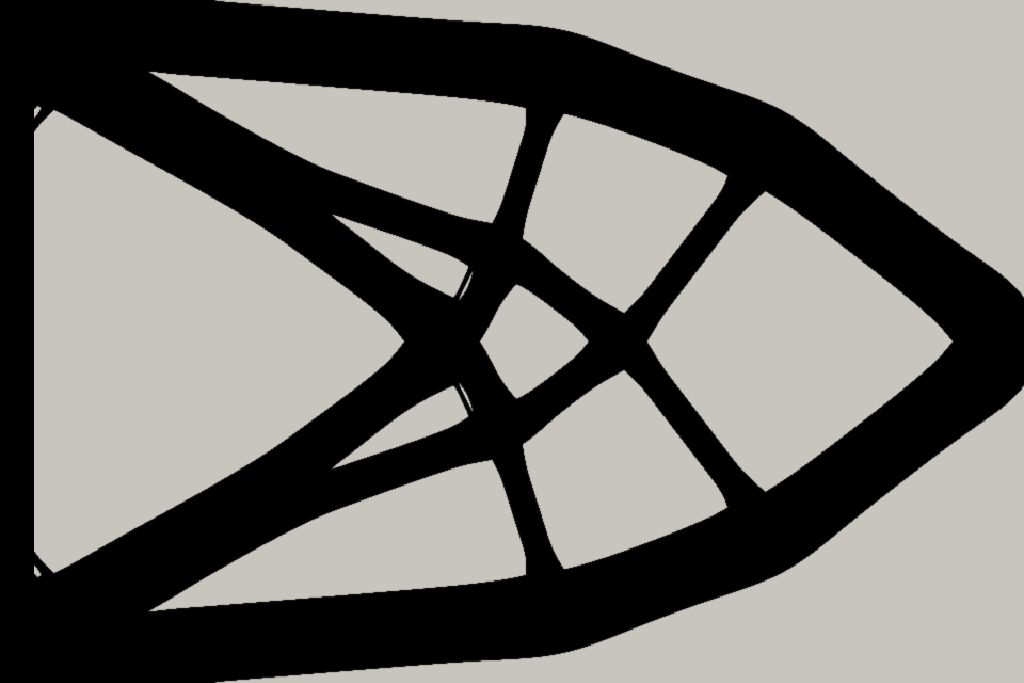}}
		 & \raisebox{-.5\height}{\includegraphics[width=0.2\textwidth]{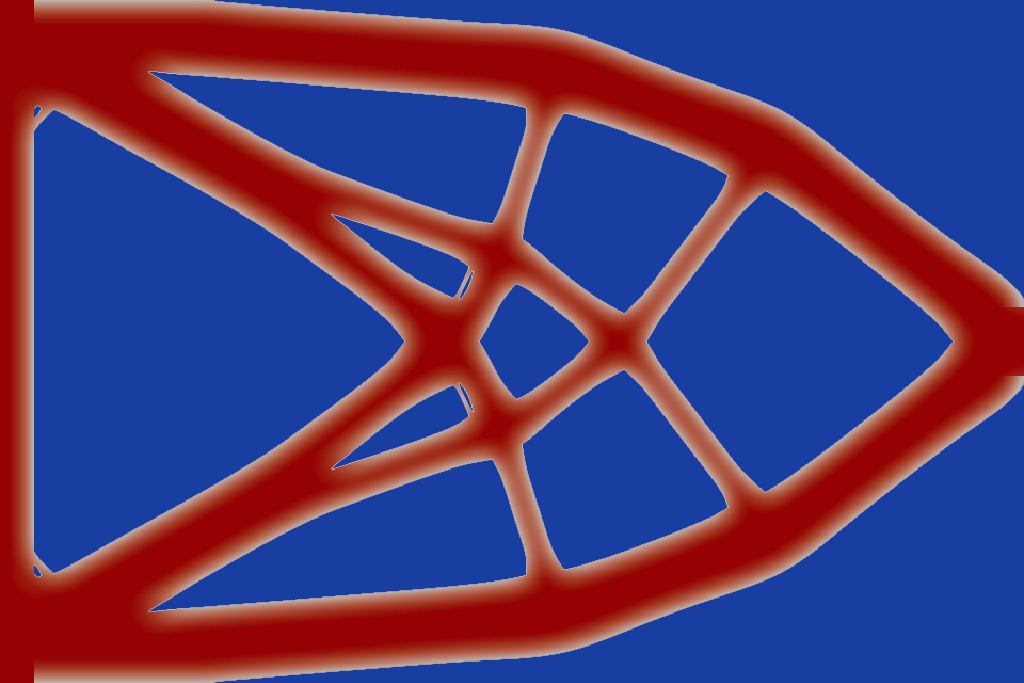}}
		\\
		 &
		 & $c=\SI{93.86}{\newton\milli\meter}$
		 &
		 & $c=\SI{85.42}{\newton\milli\meter}$
		\\[4mm]
		$\dfrac{H}{8}$
		 & \raisebox{-.5\height}{\includegraphics[width=0.2\textwidth]{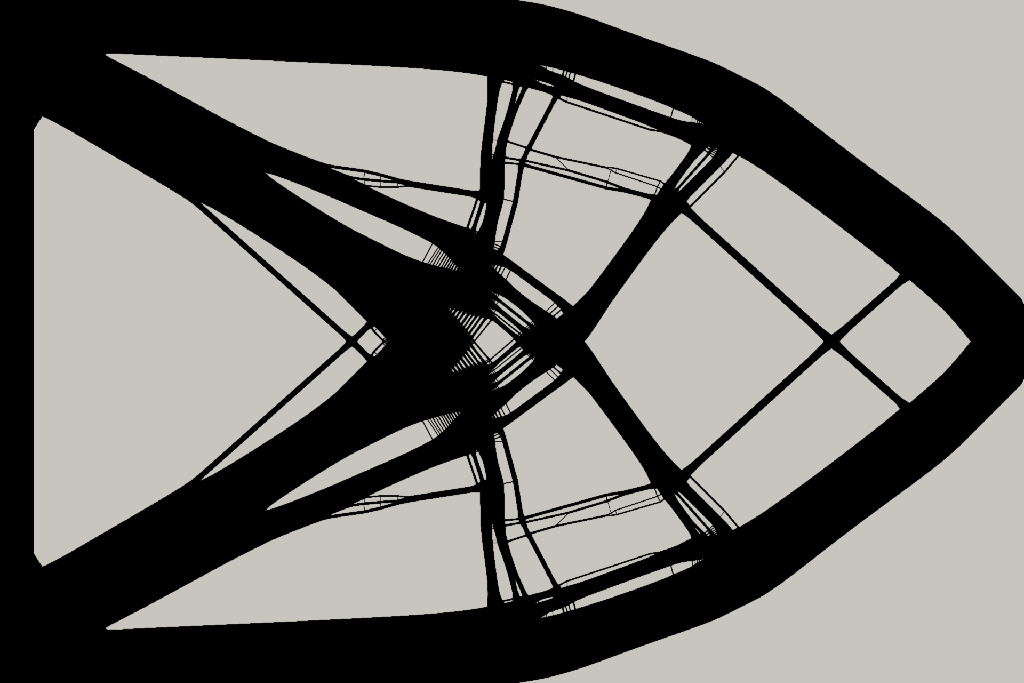}}
		 & \raisebox{-.5\height}{\includegraphics[width=0.2\textwidth]{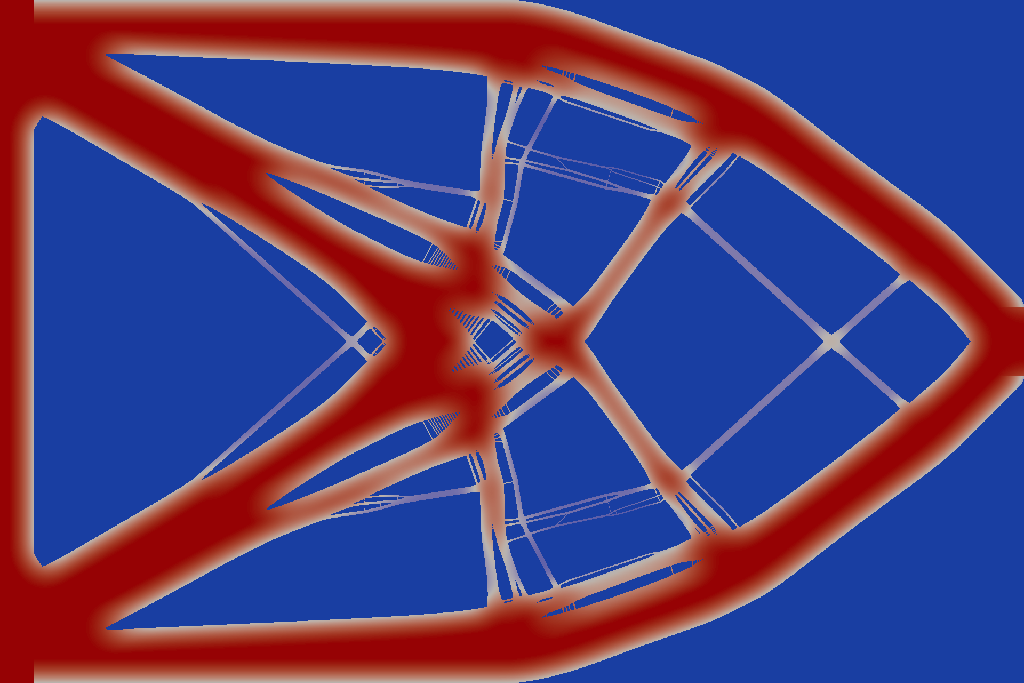}}
		 & \raisebox{-.5\height}{\includegraphics[width=0.2\textwidth]{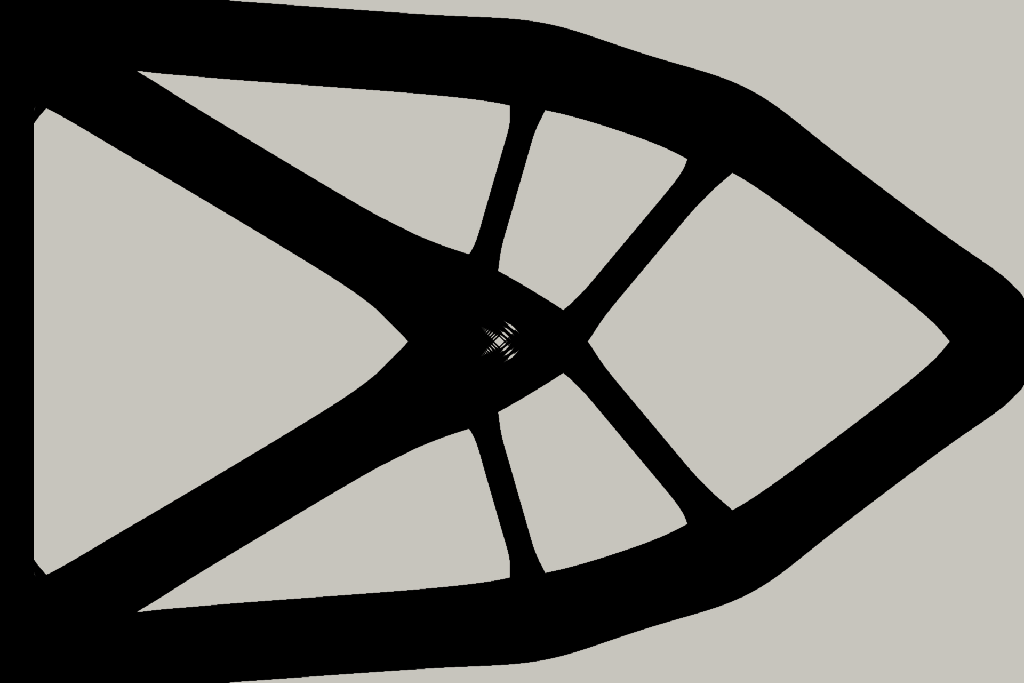}}
		 & \raisebox{-.5\height}{\includegraphics[width=0.2\textwidth]{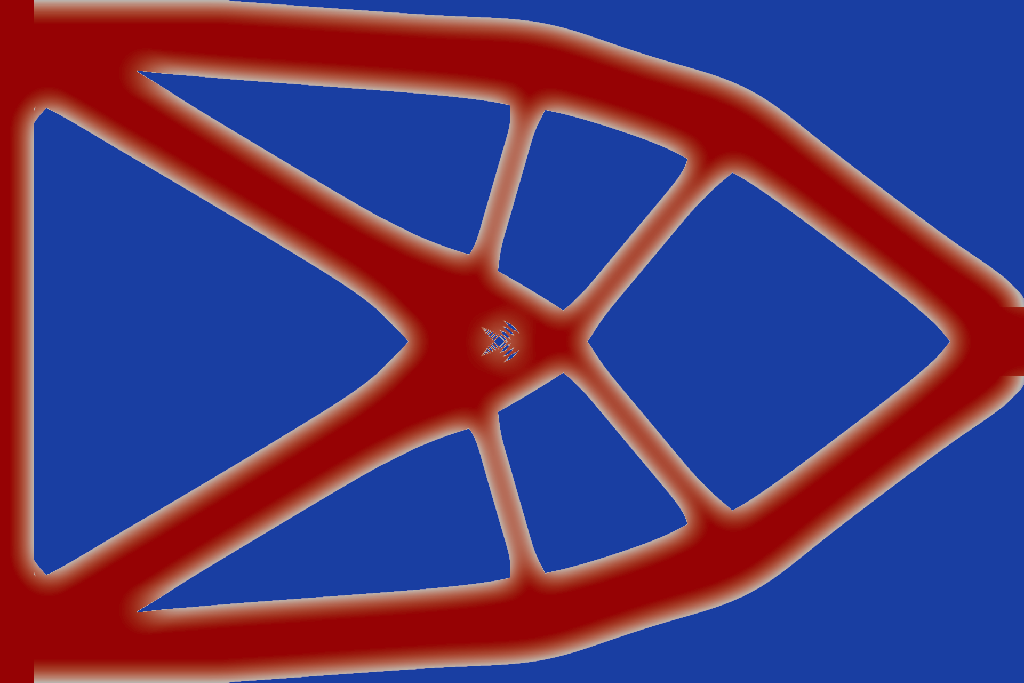}}
		\\
		 &
		 & $c=\SI{99.40}{\newton\milli\meter}$
		 &
		 & $c=\SI{87.10}{\newton\milli\meter}$
	\end{tabular}
	\vskip 0.1in
	\caption{
		Final designs optimized with $\delta=0$ and $\delta=\SI{0.4}{\milli\meter}$, and without using the Heaviside projection filter,
		for various mesh sizes $h$.
		The coarsest mesh has size $h=H=\SI{0.0625}{\milli\meter}$ and the filter radius $R$ is chosen to be $1.3h$.
		The figure shows the material distribution and the corresponding Young's modulus.
	}
	\label{fig:cantilever2d_without_heaviside}
\end{figure*}
\begin{figure*}[htbp]
	\centering
	\begin{tabular}[c]{c|cc|cc}
		\multirow{2}{*}{\begin{tabular}{c} \\ Mesh \\ size \\ ($h$) \end{tabular}}
		 & \multicolumn{2}{c}{Optimized with $\delta=0$}
		 & \multicolumn{2}{c}{Optimized with $\delta=\SI{0.4}{\milli\meter}$}
		\\[0.3cm]
		 & \begin{tabular}{c} Material \\ distribution \\ ($\rho$) \end{tabular}
		 & \begin{tabular}{c} Young's modulus \\ ($\alpha$) \end{tabular}
		 & \begin{tabular}{c} Material \\ distribution \\ ($\rho$) \end{tabular}
		 & \begin{tabular}{c} Young's modulus \\ ($\alpha$) \end{tabular}
		\\
		\hline
		\\
		$H$
		 & \raisebox{-.5\height}{\includegraphics[trim=0.5\imagewidth{} 0 0 0, clip,width=0.2\textwidth]{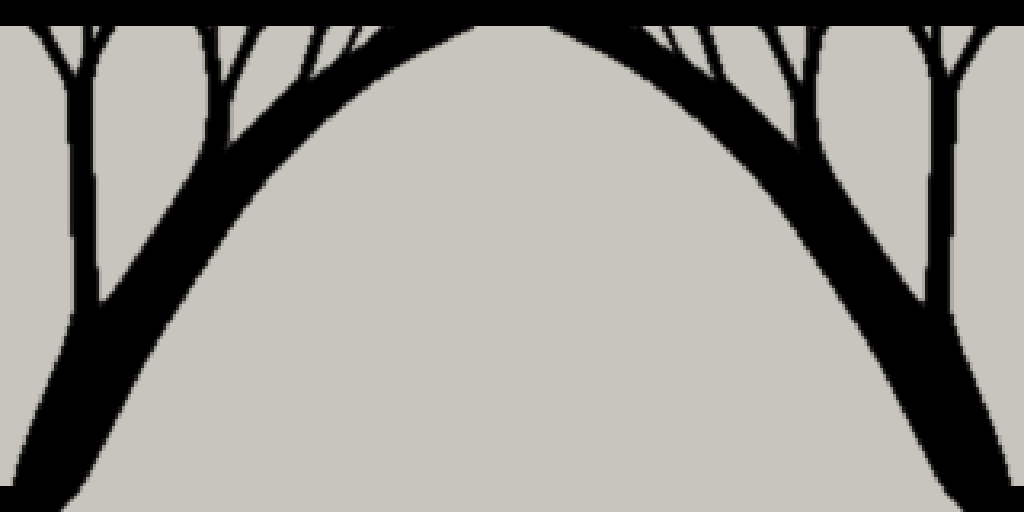}}
		 & \raisebox{-.5\height}{\includegraphics[trim=0.5\imagewidth{} 0 0 0, clip,width=0.2\textwidth]{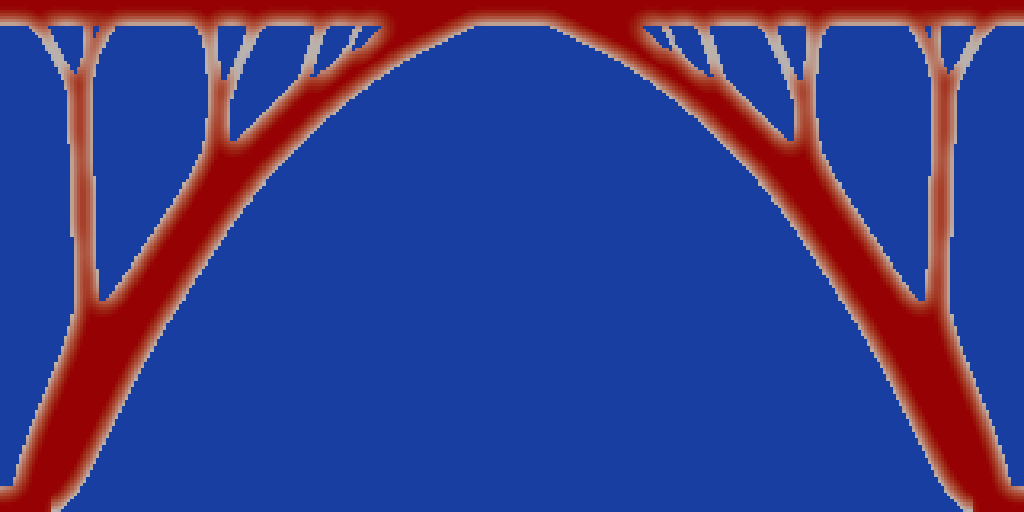}}
		 & \raisebox{-.5\height}{\includegraphics[trim=0.5\imagewidth{} 0 0 0, clip,width=0.2\textwidth]{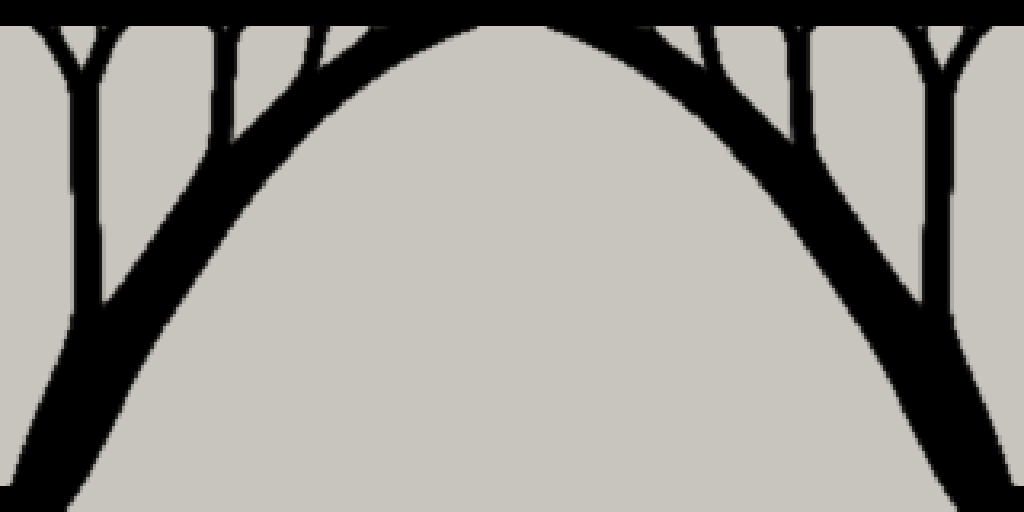}}
		 & \raisebox{-.5\height}{\includegraphics[trim=0.5\imagewidth{} 0 0 0, clip,width=0.2\textwidth]{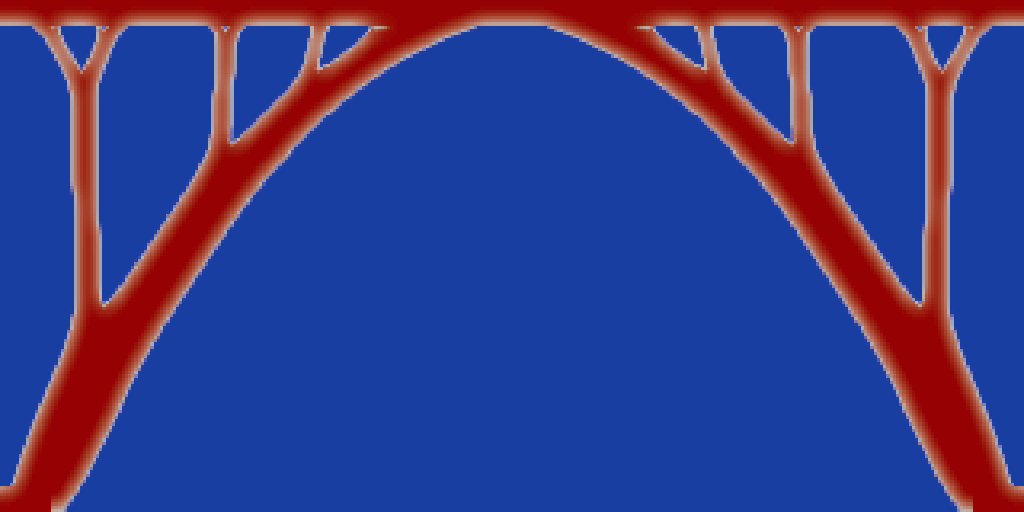}}
		\\
		 &
		 & $c=\SI{1970.77}{\newton\milli\meter}$
		 &
		 & $c=\SI{1947.64}{\newton\milli\meter}$
		\\[4mm]
		$\dfrac{H}{2}$
		 & \raisebox{-.5\height}{\includegraphics[trim=0.5\imagewidth{} 0 0 0, clip,width=0.2\textwidth]{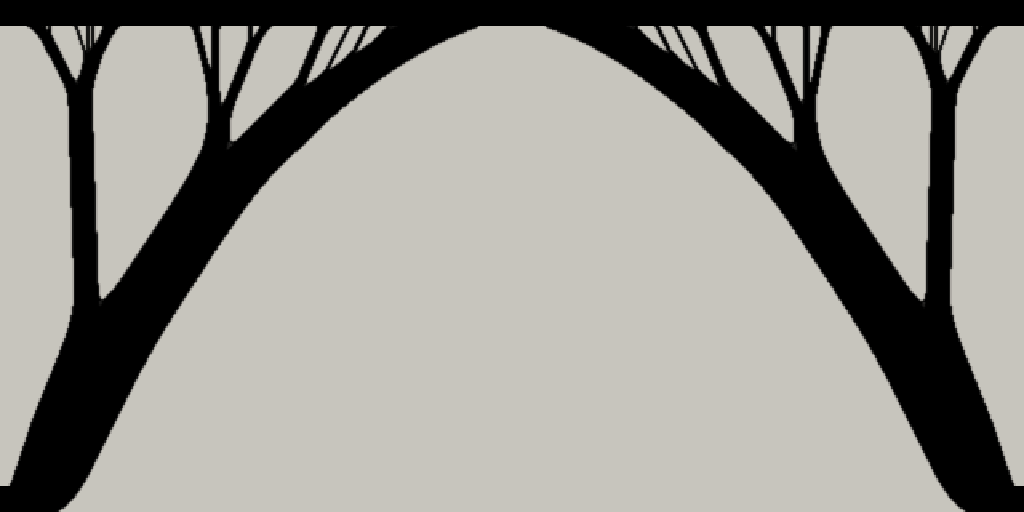}}
		 & \raisebox{-.5\height}{\includegraphics[trim=0.5\imagewidth{} 0 0 0, clip,width=0.2\textwidth]{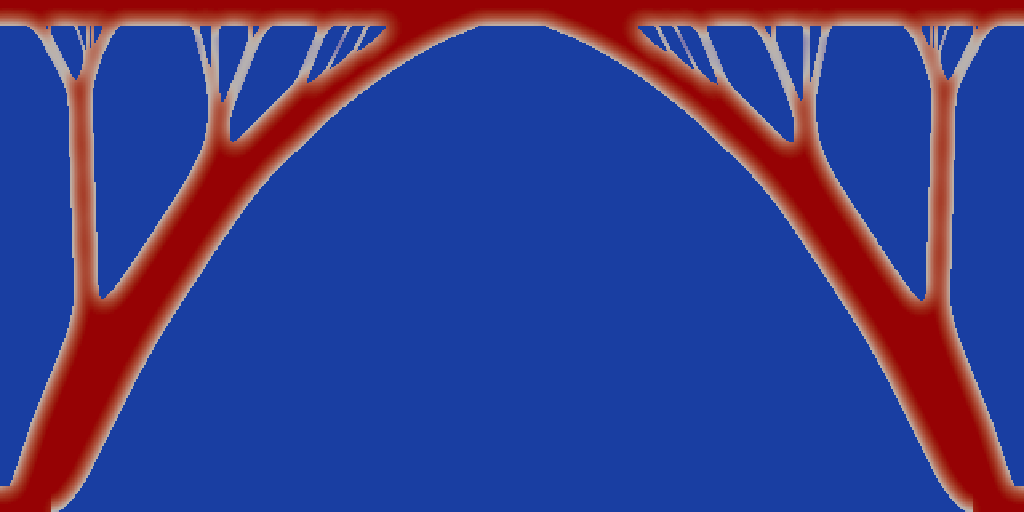}}
		 & \raisebox{-.5\height}{\includegraphics[trim=0.5\imagewidth{} 0 0 0, clip,width=0.2\textwidth]{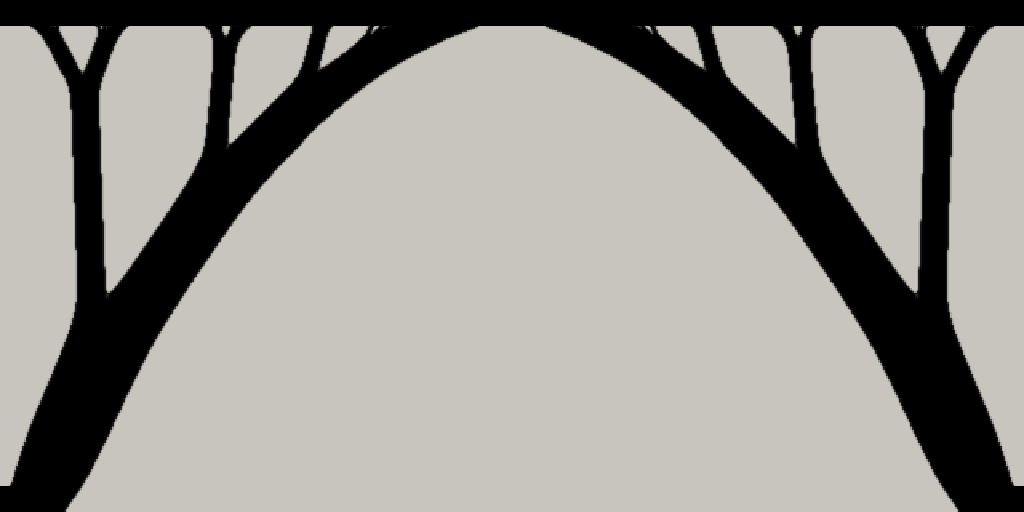}}
		 & \raisebox{-.5\height}{\includegraphics[trim=0.5\imagewidth{} 0 0 0, clip,width=0.2\textwidth]{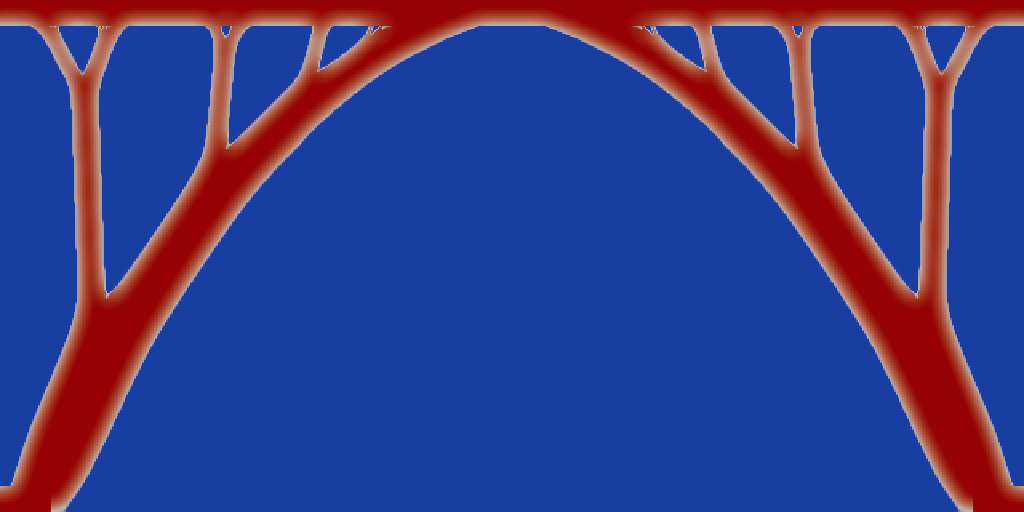}}
		\\
		 &
		 & $c=\SI{1985.92}{\newton\milli\meter}$
		 &
		 & $c=\SI{1944.03}{\newton\milli\meter}$
		\\[4mm]
		$\dfrac{H}{4}$
		 & \raisebox{-.5\height}{\includegraphics[trim=0.5\imagewidth{} 0 0 0, clip,width=0.2\textwidth]{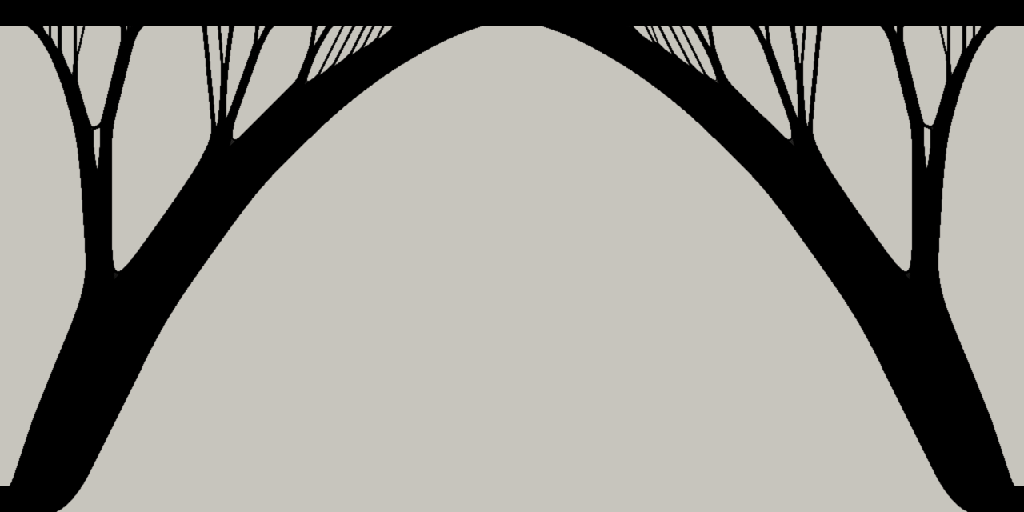}}
		 & \raisebox{-.5\height}{\includegraphics[trim=0.5\imagewidth{} 0 0 0, clip,width=0.2\textwidth]{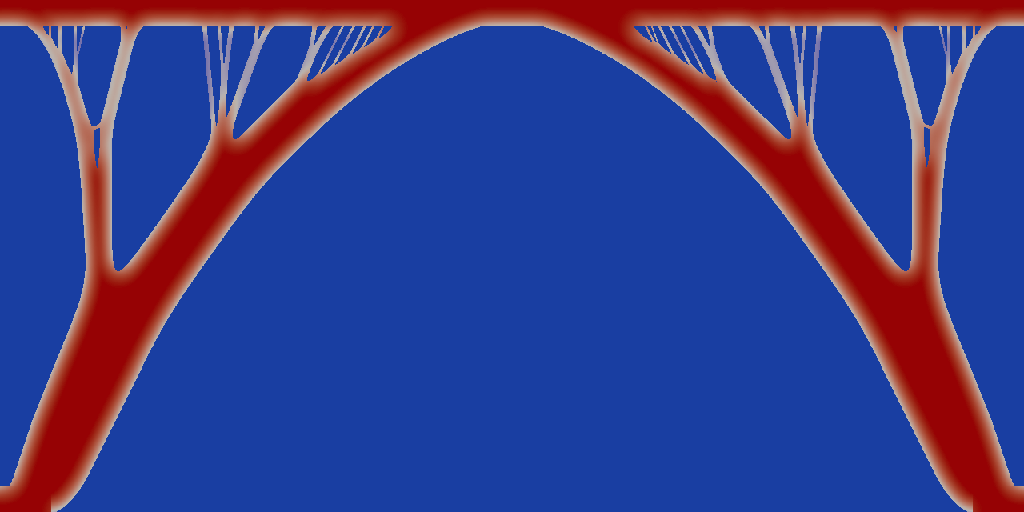}}
		 & \raisebox{-.5\height}{\includegraphics[trim=0.5\imagewidth{} 0 0 0, clip,width=0.2\textwidth]{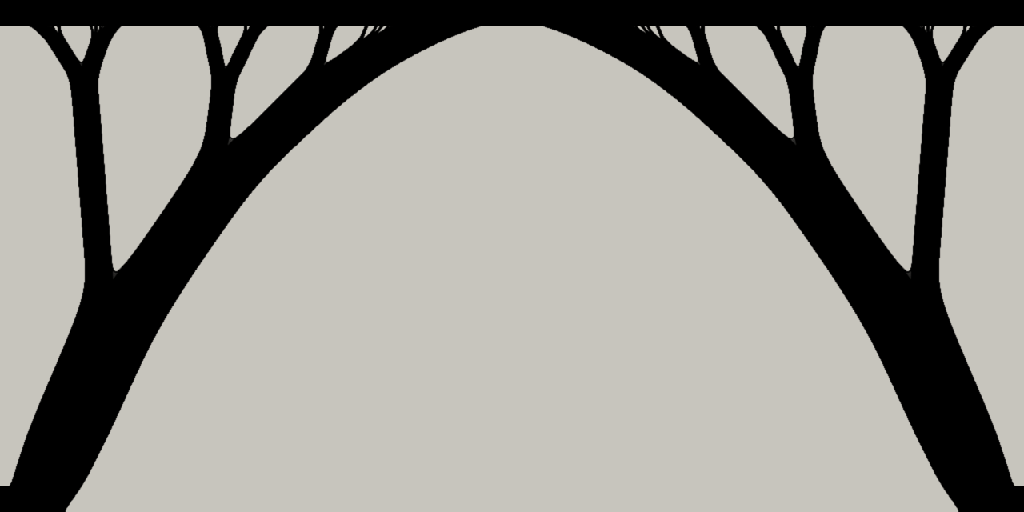}}
		 & \raisebox{-.5\height}{\includegraphics[trim=0.5\imagewidth{} 0 0 0, clip,width=0.2\textwidth]{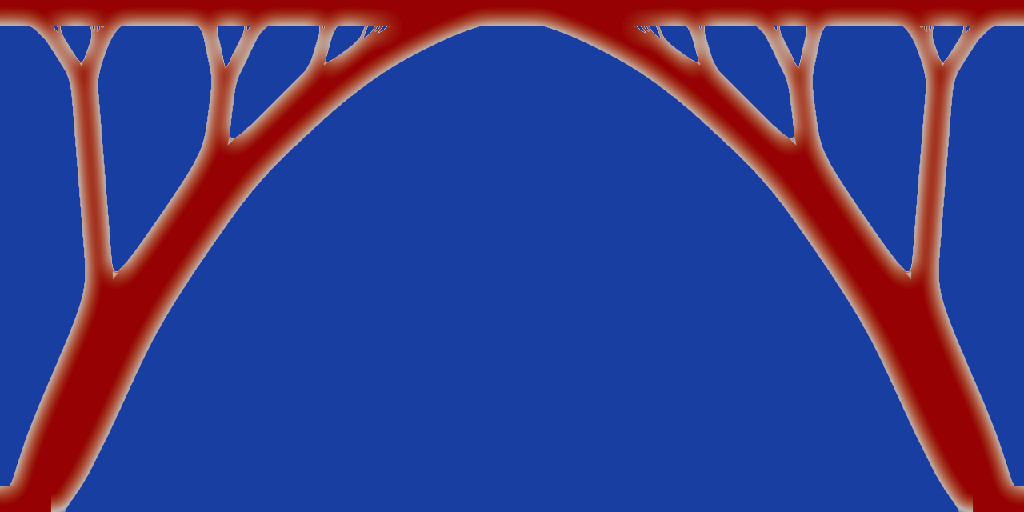}}
		\\
		 &
		 & $c=\SI{2030.03}{\newton\milli\meter}$
		 &
		 & $c=\SI{1955.86}{\newton\milli\meter}$
		\\[4mm]
		$\dfrac{H}{8}$
		 & \raisebox{-.5\height}{\includegraphics[trim=0.5\imagewidth{} 0 0 0, clip,width=0.2\textwidth]{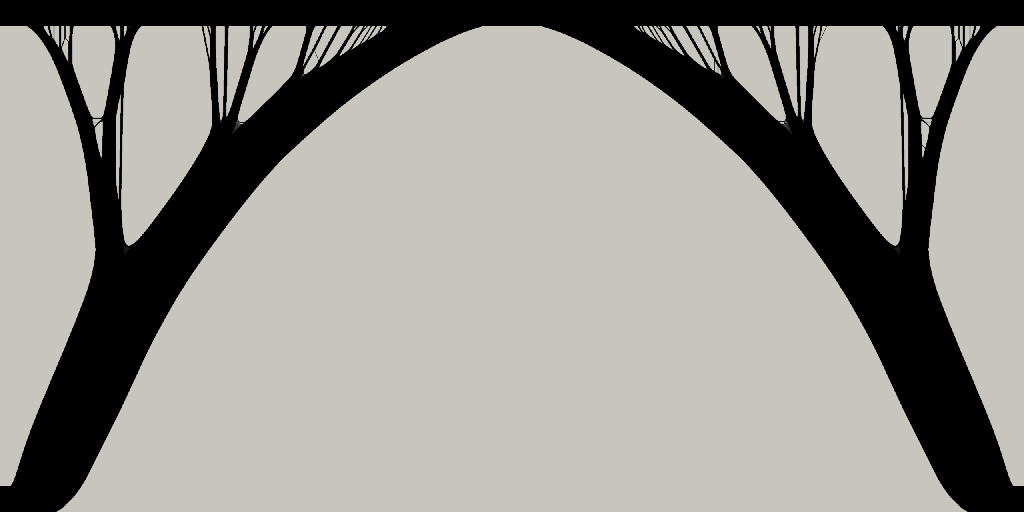}}
		 & \raisebox{-.5\height}{\includegraphics[trim=0.5\imagewidth{} 0 0 0, clip,width=0.2\textwidth]{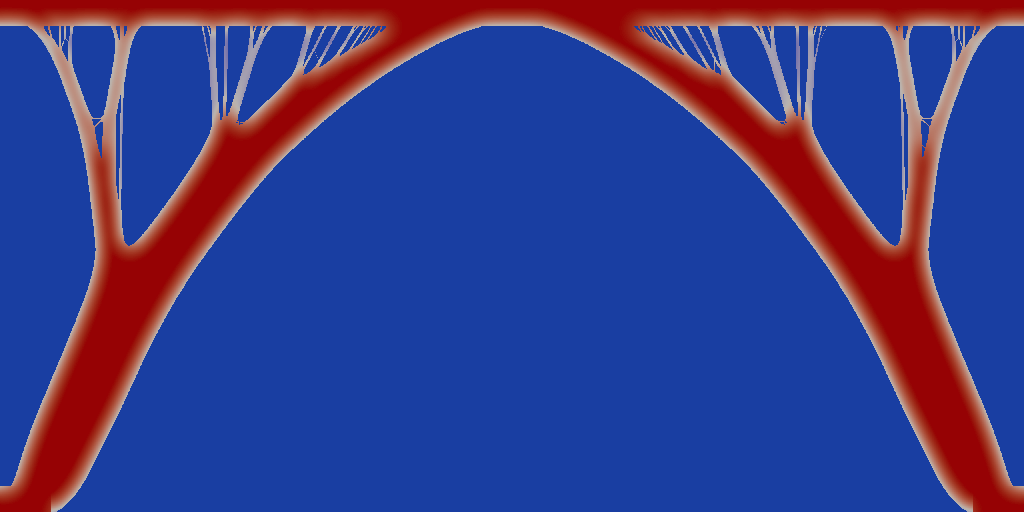}}
		 & \raisebox{-.5\height}{\includegraphics[trim=0.5\imagewidth{} 0 0 0, clip,width=0.2\textwidth]{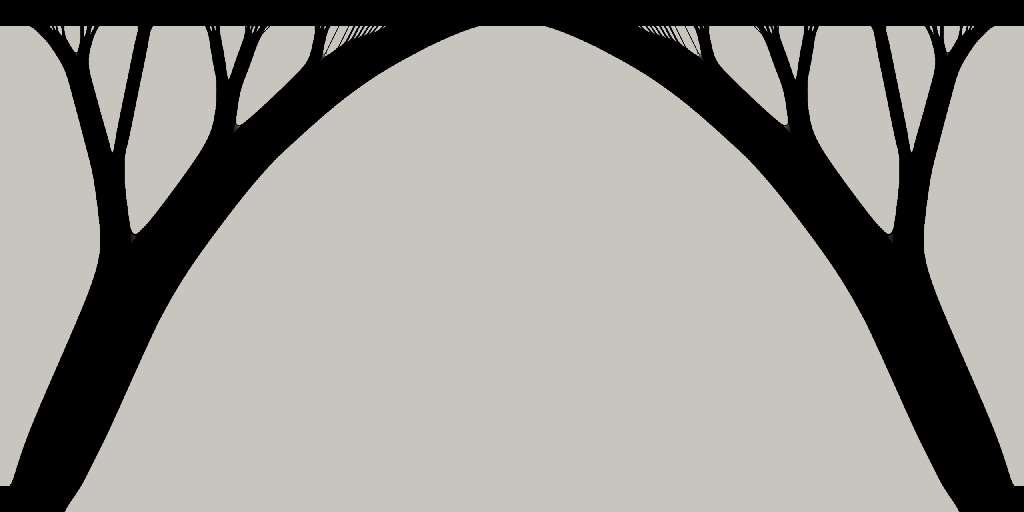}}
		 & \raisebox{-.5\height}{\includegraphics[trim=0.5\imagewidth{} 0 0 0, clip,width=0.2\textwidth]{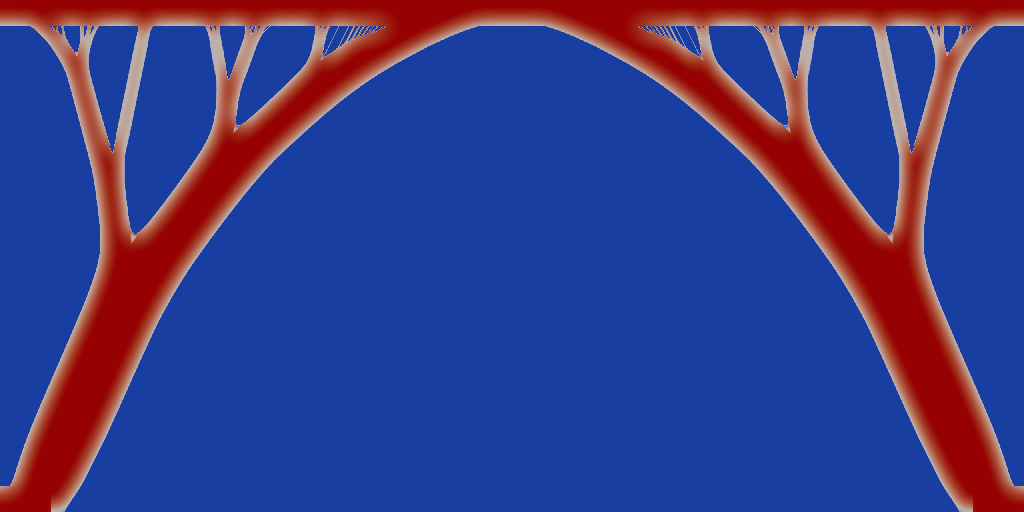}}
		\\
		 &
		 & $c=\SI{2032.8}{\newton\milli\meter}$
		 &
		 & $c=\SI{1990.98}{\newton\milli\meter}$
	\end{tabular}
	\vskip 0.1in
	\caption{
		Final designs optimized with $\delta=0$ and $\delta=\SI{0.4}{\milli\meter}$, and without using the Heaviside projection filter, for various mesh sizes $h$.
		The coarsest mesh has size $h=H=\SI{0.0625}{\milli\meter}$ and the filter radius $R$ is chosen to be $1.3h$.
		The figure shows the material distribution and the corresponding Young's modulus.
		Owing to the bridge's symmetry across the central vertical axis, only the right half of the structure is depicted.
	}
	\label{fig:bridge2d_without_heaviside}
\end{figure*}

\subsection{Mesh dependency}

Next, we study optimized designs without incorporating the Heaviside projection filter, and the necessity for a continuation scheme for $\beta$ in the optimization procedure.
While this method yields designs that exhibit sensitivity to mesh variations and lack integration of feature size control, it concurrently offers the advantage of a simplified algorithm.
By setting $R$ to a small yet sufficiently large value (e.g., $R=1.3h$), checkerboarding issues are effectively mitigated

Figure~\ref{fig:cantilever2d_without_heaviside} and Figure~\ref{fig:bridge2d_without_heaviside} visually presents the optimized designs for the cantilever beam and bridge, respectively.
These designs are generated under varying degrees of mesh refinement and are juxtaposed against the standard SIMP method, corresponding to $\delta=0$.
Numerical experiments affirm that, with moderately refined meshes, the resultant designs exhibit reduced variability when compared to those generated by the standard SIMP method.

However, there are two main downsides when using highly refined meshes. First, there's a higher chance of getting stuck in local minima, making optimization more challenging. Second, as the mesh gets finer, we start to see more complex structures (rank-2 microstructures) in the designs. This complexity adds challenges to the optimization process.

\section{Conclusion}
\label{sec:conclusion}

In this work, we have presented a framework for the topology optimization of structures featuring graded surfaces.
Our approach relies on a nonlocal methodology that incorporates surface grading, addressing a critical aspect of additively manufactured structures.
Through the implementation of this approach, we have explored the distribution of Young's modulus across the thickness of one- and two-dimensional specimens, unveiling nuanced variations in material properties.

The primary advantage of our proposed framework is its ability to optimize structures with graded surfaces, reflecting the actual material properties of additively manufactured structures. This capability is not readily achievable with conventional topology optimization methods using the SIMP formulation.
Moreover, our approach, particularly in its simplified form, eliminates the need for a continuation scheme, which is often employed in SIMP-based optimization schemes using Heaviside projection.

In summary, the findings suggest that incorporating actual material information into the optimization process leads to better designs in terms of compliance, particularly in smaller-scale, thin-walled structures.
While the difference in compliance improvement may not be substantial in absolute terms, the trend towards enhanced stiffness and reduced deformation with topology optimization is consistent with theoretical expectations and offers a valuable avenue for optimizing structural performance.

The present study can be extended in several directions, for example, adding stress and buckling constraints to the optimization problem.
Additionally, the introduction of material anisotropy, as explored in \cite{TopologyOptimiSuresh2020}, could further enhance the versatility and applicability of the proposed framework.

\subsection*{Replication of results}
To facilitate the replication of the results presented in this paper, the accompanying code is publicly available on GitHub at \url{https://github.com/sukhmindersingh/topopt_nonlocal}.

\subsection*{Competing interests and funding}
This study was funded by the Deutsche Forschungsgemeinschaft (DFG, German Research Foundation) -- Project-ID 61375930 -- SFB 814, sub-project T3, A3.
The authors have no competing interests to declare that are relevant to the content of this article.

\bibliography{references}
\end{document}

%% file: packages.tex
\usepackage[utf8]{inputenc}
\usepackage[T1]{fontenc}
\usepackage[english]{babel}

\usepackage{lineno}
\usepackage{standalone}
\usepackage{caption}
\usepackage{subcaption}
\usepackage{stmaryrd}
\usepackage{amssymb,amsmath,amstext,empheq}
\usepackage{upgreek}
\usepackage{graphicx}
\usepackage{booktabs}
\usepackage{float}

\usepackage{tikz}
\usetikzlibrary{decorations.pathmorphing,decorations.pathreplacing,decorations.markings,patterns,calc}
\usetikzlibrary{shapes.geometric, shapes.arrows, arrows, angles, quotes}
\usetikzlibrary{spy,backgrounds}
\usetikzlibrary{external}
\usetikzlibrary{intersections}
\usetikzlibrary{tikzmark}

\tikzstyle{startstop} = [rectangle, rounded corners, minimum width=3cm, minimum height=1cm,text centered, draw=black, fill=red!30]
\tikzstyle{io} = [trapezium, trapezium left angle=70, trapezium right angle=110, minimum width=3cm, minimum height=1cm, text centered, draw=black, fill=blue!30]
\tikzstyle{process} = [rectangle, minimum width=3cm, minimum height=1cm, text centered, draw=black, fill=orange!30]
\tikzstyle{decision} = [diamond, minimum width=3cm, minimum height=1cm, text centered, draw=black, fill=green!30]
\tikzstyle{arrow} = [->,>=stealth]

\usepackage{import}
\usepackage{pdfpages}
\usepackage{transparent}
\usepackage{xcolor}

\usepackage[normalem]{ulem}
\usepackage{textcomp}
\usepackage{siunitx}
\usepackage[most]{tcolorbox}

\usepackage{bm}
\usepackage[]{datatool}
\usepackage{calc}
\usepackage{textgreek}
\usepackage{bigints}
\usepackage{physics}
\usepackage{interval}
\usepackage{accents}
\usepackage{mathtools,eqparbox}
\usepackage{indentfirst}
\usepackage[algo2e,ruled,linesnumbered,vlined,commentsnumbered,longend]{algorithm2e}
\SetKw{Continue}{continue}
\SetKw{Break}{break}

\usepackage{multirow}
\usepackage{rotating}
\usepackage{adjustbox}
\usepackage{array}
\usepackage{enumerate}

\usepackage{pgfplots}
\pgfplotsset{compat=newest}
\usepgfplotslibrary{groupplots,dateplot}
\usepgfplotslibrary{colormaps}
\usepgfplotslibrary{colorbrewer}
\usepgfplotslibrary{polar}
\usepgfplotslibrary{fillbetween}

\usepackage[absolute,overlay]{textpos}

\usepackage{tikz-dimline}

\usepackage[autostyle]{csquotes}

\usepackage[colorinlistoftodos]{todonotes}

\usepackage{hyperref}
\usepackage{hypcap}
\usepackage{bookmark}

\newcommand{\rom}[1]{\uppercase\expandafter{\romannumeral #1\relax}}

\intervalconfig {
	soft open fences ,
}

\usepackage{amsthm}

\theoremstyle{plain}

\theoremstyle{definition}

\theoremstyle{remark}

\pgfplotsset{
	legend image with text/.style={
			legend image code/.code={%
					\node[anchor=center] at (0.3cm,0cm) {#1};
				}
		},
}

%% file: main.bbl
\begin{thebibliography}{}
\providecommand{\doi}[1]{\url{https://doi.org/#1}}
\bibcommenthead

\bibitem[\protect\citeauthoryear{Algardh, Horn, West, Aman, Snis, Engqvist, Lausmaa, and Harrysson}{Algardh et~al.}{2016}]{ThicknessDepenAlgard2016}
Algardh, J.K., T.~Horn, H.~West, R.~Aman, A.~Snis, H.~Engqvist, J.~Lausmaa, and O.~Harrysson. 2016, 10.
\newblock Thickness dependency of mechanical properties for thin-walled titanium parts manufactured by {Electron Beam Melting} ({EBM}).
\newblock {\em Additive Manufacturing\/}~12: 45--50.
\newblock \doi{10.1016/j.addma.2016.06.009} .

\bibitem[\protect\citeauthoryear{Bendsøe and Sigmund}{Bendsøe and Sigmund}{2004}]{TopologyOptimiBendso2004}
Bendsøe, M.P. and O.~Sigmund. 2004.
\newblock {\em Topology Optimization}.
\newblock Springer Berlin Heidelberg.

\bibitem[\protect\citeauthoryear{Bruns and Tortorelli}{Bruns and Tortorelli}{1998}]{TopologyOptimiBruns1998}
Bruns, T. and D.~Tortorelli 1998, 9.
\newblock Topology optimization of geometrically nonlinear structures and compliant mechanisms.
\newblock In {\em 7th AIAA/USAF/NASA/ISSMO Symposium on Multidisciplinary Analysis and Optimization}. American Institute of Aeronautics and Astronautics.

\bibitem[\protect\citeauthoryear{Chin, Leader, and Kennedy}{Chin et~al.}{2019}]{AScalableFramChin2019}
Chin, T.W., M.K. Leader, and G.J. Kennedy. 2019, 9.
\newblock A scalable framework for large-scale 3d multimaterial topology optimization with octree-based mesh adaptation.
\newblock {\em Advances in Engineering Software\/}~135: 102682.
\newblock \doi{10.1016/j.advengsoft.2019.05.004} .

\bibitem[\protect\citeauthoryear{Clausen, Aage, and Sigmund}{Clausen et~al.}{2015}]{TopologyOptimiClause2015}
Clausen, A., N.~Aage, and O.~Sigmund. 2015, 6.
\newblock Topology optimization of coated structures and material interface problems.
\newblock {\em Computer Methods in Applied Mechanics and Engineering\/}~290: 524--541.
\newblock \doi{10.1016/j.cma.2015.02.011} .

\bibitem[\protect\citeauthoryear{Eringen}{Eringen}{1987}]{TheoryOfNonloEringe1987}
Eringen, A.C. 1987.
\newblock Theory of nonlocal elasticity and some applications.
\newblock {\em Res Mech.\/}~{\em 21\/}(4): 313--342 .

\bibitem[\protect\citeauthoryear{Groen, Wu, and Sigmund}{Groen et~al.}{2019}]{HomogenizationGroen2019}
Groen, J.P., J.~Wu, and O.~Sigmund. 2019, 6.
\newblock Homogenization-based stiffness optimization and projection of 2d coated structures with orthotropic infill.
\newblock {\em Computer Methods in Applied Mechanics and Engineering\/}~349: 722--742.
\newblock \doi{10.1016/j.cma.2019.02.031} .

\bibitem[\protect\citeauthoryear{Guest, Prévost, and Belytschko}{Guest et~al.}{2004}]{AchievingMinimGuest2004}
Guest, J.K., J.H. Prévost, and T.~Belytschko. 2004, 9.
\newblock Achieving minimum length scale in topology optimization using nodal design variables and projection functions.
\newblock {\em International Journal for Numerical Methods in Engineering\/}~61: 238--254.
\newblock \doi{10.1002/nme.1064} .

\bibitem[\protect\citeauthoryear{Jaksch, Spinola, Cholewa, Pflug, Stingl, and Drummer}{Jaksch et~al.}{2022}]{jaksch2022thin}
Jaksch, A., M.~Spinola, C.~Cholewa, L.~Pflug, M.~Stingl, and D.~Drummer 2022.
\newblock Thin-walled part properties in {PBF-LB/P}—experimental understanding and nonlocal material model.
\newblock In {\em 2022 International Solid Freeform Fabrication Symposium}.

\bibitem[\protect\citeauthoryear{Luo, Li, and Liu}{Luo et~al.}{2019}]{AProjectionBaLuoY2019}
Luo, Y., Q.~Li, and S.~Liu. 2019, 6.
\newblock A projection‐based method for topology optimization of structures with graded surfaces.
\newblock {\em International Journal for Numerical Methods in Engineering\/}~118: 654--677.
\newblock \doi{10.1002/nme.6031} .

\bibitem[\protect\citeauthoryear{Sigmund}{Sigmund}{2007}]{MorphologyBaseSigmun2007}
Sigmund, O. 2007, 2.
\newblock Morphology-based black and white filters for topology optimization.
\newblock {\em Structural and Multidisciplinary Optimization\/}~33: 401--424.
\newblock \doi{10.1007/s00158-006-0087-x} .

\bibitem[\protect\citeauthoryear{Sindinger, Kralovec, Tasch, and Schagerl}{Sindinger et~al.}{2020}]{ThicknessDepenSindin2020}
Sindinger, S.L., C.~Kralovec, D.~Tasch, and M.~Schagerl. 2020, 5.
\newblock Thickness dependent anisotropy of mechanical properties and inhomogeneous porosity characteristics in laser-sintered polyamide 12 specimens.
\newblock {\em Additive Manufacturing\/}~33: 101141.
\newblock \doi{10.1016/j.addma.2020.101141} .

\bibitem[\protect\citeauthoryear{Sindinger, Marschall, Kralovec, and Schagerl}{Sindinger et~al.}{2021a}]{ConsideringInhSindin2021}
Sindinger, S.L., D.~Marschall, C.~Kralovec, and M.~Schagerl. 2021a.
\newblock Considering inhomogeneous material properties for stiffness and failure prediction of thin-walled additively manufactured parts.
\newblock {\em Procedia Structural Integrity\/}~34: 78--86.
\newblock \doi{10.1016/j.prostr.2021.12.012} .

\bibitem[\protect\citeauthoryear{Sindinger, Marschall, Kralovec, and Schagerl}{Sindinger et~al.}{2021b}]{MaterialModellSindin2021}
Sindinger, S.L., D.~Marschall, C.~Kralovec, and M.~Schagerl. 2021b, 1.
\newblock Material modelling and property mapping for structural fea of thin-walled additively manufactured components.
\newblock {\em Virtual and Physical Prototyping\/}~16: 97--112.
\newblock \doi{10.1080/17452759.2020.1824427} .

\bibitem[\protect\citeauthoryear{Sindinger, Marschall, Kralovec, and Schagerl}{Sindinger et~al.}{2021c}]{StructuralRespSindin2021}
Sindinger, S.L., D.~Marschall, C.~Kralovec, and M.~Schagerl. 2021c, 5.
\newblock Structural response prediction of thin-walled additively manufactured parts considering orthotropy, thickness dependency and scatter.
\newblock {\em Materials\/}~14: 2463.
\newblock \doi{10.3390/ma14092463} .

\bibitem[\protect\citeauthoryear{Suresh, Thore, Torstenfelt, and Klarbring}{Suresh et~al.}{2020}]{TopologyOptimiSuresh2020}
Suresh, S., C.J. Thore, B.~Torstenfelt, and A.~Klarbring. 2020, 12.
\newblock Topology optimization accounting for surface layer effects.
\newblock {\em Structural and Multidisciplinary Optimization\/}~62: 3009--3019.
\newblock \doi{10.1007/s00158-020-02644-x} .

\bibitem[\protect\citeauthoryear{Svanberg}{Svanberg}{1987}]{TheMethodOfMSvanbe1987}
Svanberg, K. 1987, 2.
\newblock The method of moving asymptotes—a new method for structural optimization.
\newblock {\em International Journal for Numerical Methods in Engineering\/}~24: 359--373.
\newblock \doi{10.1002/nme.1620240207} .

\bibitem[\protect\citeauthoryear{Tasch, Mad, Stadlbauer, and Schagerl}{Tasch et~al.}{2018}]{ThicknessDepenTasch2018}
Tasch, D., A.~Mad, R.~Stadlbauer, and M.~Schagerl. 2018, 10.
\newblock Thickness dependency of mechanical properties of laser-sintered polyamide lightweight structures.
\newblock {\em Additive Manufacturing\/}~23: 25--33.
\newblock \doi{10.1016/j.addma.2018.06.018} .

\bibitem[\protect\citeauthoryear{Tuna and Trovalusci}{Tuna and Trovalusci}{2022}]{TopologyOptimiTuna2022}
Tuna, M. and P.~Trovalusci. 2022, 9.
\newblock Topology optimization of scale-dependent non-local plates.
\newblock {\em Structural and Multidisciplinary Optimization\/}~65.
\newblock \doi{10.1007/s00158-022-03351-5} .

\bibitem[\protect\citeauthoryear{Wang and Kang}{Wang and Kang}{2018}]{ALevelSetMetWang2018}
Wang, Y. and Z.~Kang. 2018, 2.
\newblock A level set method for shape and topology optimization of coated structures.
\newblock {\em Computer Methods in Applied Mechanics and Engineering\/}~329: 553--574.
\newblock \doi{10.1016/j.cma.2017.09.017} .

\bibitem[\protect\citeauthoryear{Wu, Sigmund, and Groen}{Wu et~al.}{2021}]{TopologyOptimiWuJu2021}
Wu, J., O.~Sigmund, and J.P. Groen. 2021, 3.
\newblock Topology optimization of multi-scale structures: a review.
\newblock {\em Structural and Multidisciplinary Optimization\/}~63: 1455--1480.
\newblock \doi{10.1007/s00158-021-02881-8} .

\bibitem[\protect\citeauthoryear{Yoon and Yi}{Yoon and Yi}{2019}]{ANewCoatingFYoon2019}
Yoon, G.H. and B.~Yi. 2019, 10.
\newblock A new coating filter of coated structure for topology optimization.
\newblock {\em Structural and Multidisciplinary Optimization\/}~60: 1527--1544.
\newblock \doi{10.1007/s00158-019-02279-7} .

\end{thebibliography}
